%% file: main.tex
\documentclass[AMA,Times1COL]{WileyNJDv5_arxiv} 

\articletype{}%

\startpage{1}

\graphicspath{ {./figs/} }
\usepackage{multirow}
\usepackage{booktabs}
\usepackage[section]{placeins}
\usepackage{caption}

\usepackage{csml-wiley}		

\usepackage{algpseudocode}
\usepackage{algorithm}

\hypersetup{
    colorlinks=true,
    linkcolor=black,
    citecolor=black,
    filecolor=magenta,      
    urlcolor=cyan
}

\begin{document}

\title{Bayesian buckling load optimisation for structures with geometric uncertainties}

\author[1]{Tianyi Liu}

\author[2]{Xiao Xiao}

\author[1]{Fehmi Cirak}

\authormark{LIU \textsc{et al.}}
\titlemark{Bayesian buckling load optimisation for structures with geometric uncertainties}

\address[1]{\orgdiv{Department of Engineering}, \orgname{University of Cambridge}, \orgaddress{\state{Cambridge}, \country{United Kingdom}}}

\address[2]{\orgdiv{School of Ocean and Civil Engineering}, \orgname{Shanghai Jiao Tong University}, \orgaddress{\state{Shanghai}, \country{China}}}

\corres{Corresponding author Xiao Xiao, School of Ocean and Civil Engineering, Shanghai Jiao Tong University, Shanghai, China, 200240. \email{xiao.xiao@sjtu.edu.cn}}



\abstract[Abstract]
{Optimised lightweight structures, such as shallow domes and slender towers, are prone to sudden buckling failure because geometric uncertainties/imperfections can lead to a drastic reduction in their buckling loads. We introduce a framework for the robust optimisation of buckling loads, considering geometric nonlinearities and random geometric imperfections. The mean and standard deviation of buckling loads are estimated by Monte Carlo sampling of random imperfections and performing a nonlinear finite element computation for each sample. The extended system method is employed to compute the buckling load directly, avoiding costly path-following procedures. Furthermore, the quasi-Monte Carlo sampling using the Sobol sequence is implemented to generate more uniformly distributed samples, which significantly reduces the number of finite element computations. The objective function consisting of the weighted sum of the mean and standard deviation of the buckling load is optimised using Bayesian optimisation. The accuracy and efficiency of the proposed framework are demonstrated through robust sizing optimisation of several geometrically nonlinear truss examples.}


\jnlcitation{\cname{%
\author{Liu T.},
\author{Xiao X.}, and
\author{Cirak F.}}.
\ctitle{Bayesian robust optimisation of stability of truss structures with geometric imperfections.} \cjournal{\it J Comput Phys.} \cvol{2021;00(00):1--18}.}

\maketitle

\input{introduction}

\input{geoNonlinear}

\input{directRandStability}

\input{bayesianOpt}

\input{examples}

\input{conclusions}



\bmsection*{Conflict of interest}

The authors declare no potential conflict of interests.

\bibliography{references}



\end{document}

%% file: introduction.tex
\section{Introduction \label{sec:introduction}}
%
Structural optimisation has become a commonly used tool in industrial design. With the conventional deterministic optimisation of linear structures, the safety and robustness of an optimised structure are not automatically guaranteed. First and foremost, many optimised lightweight structures tend to have slender bars and thin walls, which are susceptible to sudden failure through buckling. The optimisation of linear structures does not consider structural instability and consequent failure. Hence, it is essential to consider nonlinear structural behaviour and instability in the optimisation process~\cite{Rahmatalla2003, Lindgaard2010, Madah2017}. 
Moreover, the as-built structure inevitably has random imperfections deviating from the ideal design, which are typically induced by manufacturing, assembly and construction~\cite{Thompson1969}. The design determined by conventional deterministic optimisation is usually sensitive to those uncertainties, resulting in the loss of optimality of the as-built structure. An imperfect structure often has a significantly altered equilibrium path and a greatly reduced buckling load compared to that of the perfect structure. Therefore, it is essential to include uncertainties, such as variations in geometry~\cite{Kang2013, Jansen2015, Madah2019}, material properties~\cite{Rostami2021, Bai2021, YelunCirak2024} and loading~\cite{Dunning2011, Zhao2014, Gao2022}, and perform robust structural optimisation~\cite{Schueller2001, Schueller2008, Beyer2007} to obtain a more robust design.

\begin{figure*}[hbtp!]
	\centering
	\subfloat[Top view of the design] {
		\includegraphics[width=0.23 \textwidth]{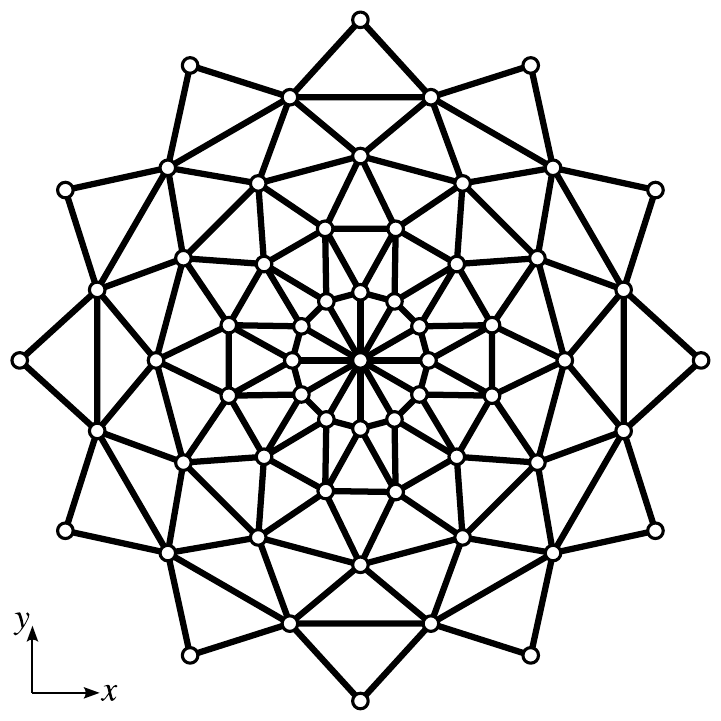}
		\label{fig:starDome5_perfect} }
	\hfil
	\subfloat[Side view of the design] {
		\includegraphics[width=0.23 \textwidth]{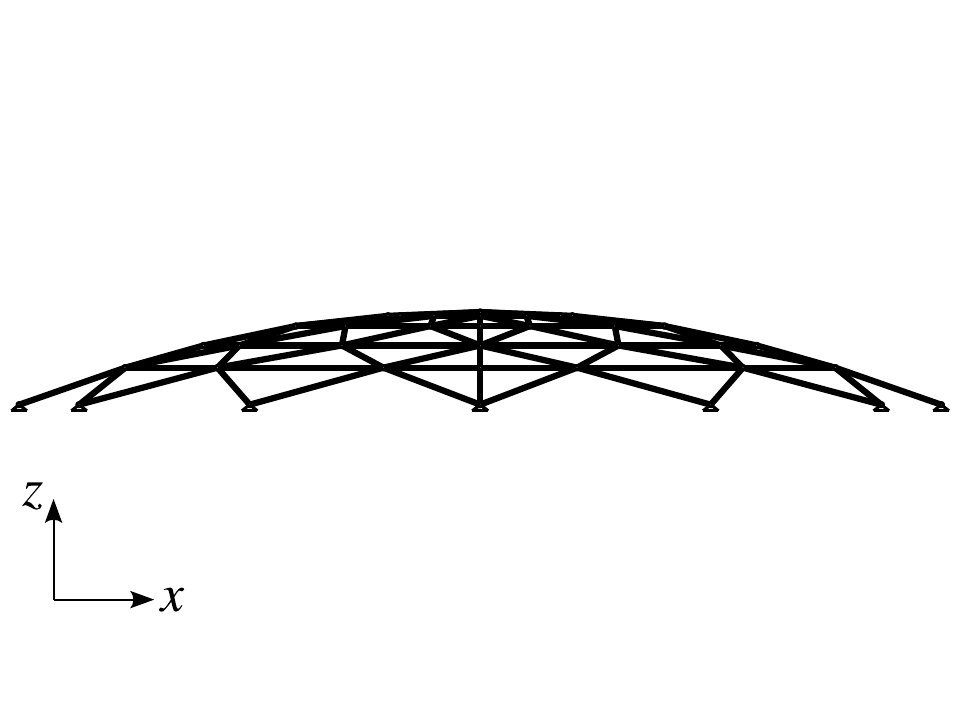}
		\label{fig:starDome5n_imperfect} }
	\hfil
	\subfloat[Mean-value optimised structure] {
		\includegraphics[width=0.23 \textwidth]{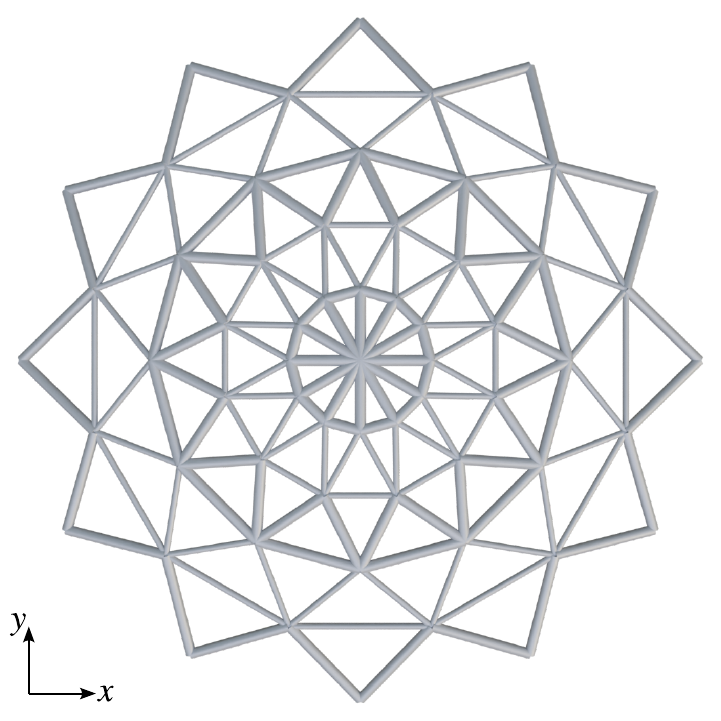}
		\label{fig:starDome5n_deterministicOptResult} }
	\hfil
	\subfloat[Mean and std optimised structure] {
		\includegraphics[width=0.23 \textwidth]{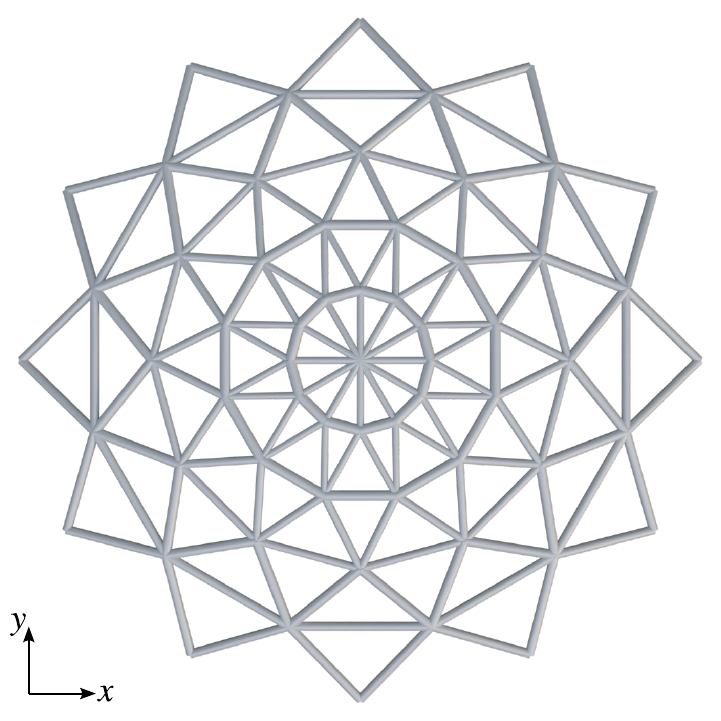}
		\label{fig:starDome5_robustOptResult} }
	\caption{Optimisation of a pin-jointed star dome truss. The aim is to maximise the mean and to minimise the standard deviation (std) of the buckling load, and the design variables are the cross-sectional areas of struts. The optimisation with respect to only the mean buckling load in figure~(c) yields the mean buckling load $\overline{\lambda}_c/\overline{\lambda}_c^* = 1.0$ and standard deviation $\sigma_{c}/\sigma_{c}^* = 1.0$. The optimisation with respect to the mean and standard deviation of the buckling load in figure~(d) yields the mean buckling load $\overline{\lambda}_c/\overline{\lambda}_c^* = 0.646$ and standard deviation $\sigma_{c}/\sigma_{c}^* = 0.491$. The two parameters~$\lambda_c^*$ and~$\sigma_{c}^*$ are normalising constants. For further details see Section~\ref{sec:examples}.}
	\label{fig:starDome5_imperfect_robustOpt}
\end{figure*}

In this paper, we propose a framework for robust optimisation of buckling loads by taking into account nonlinearities and geometric uncertainties. The framework is applied to robust sizing optimisation of pin-jointed spatial trusses with geometric nonlinearities and geometric imperfections, which are random; see Figure~\ref{fig:starDome5_imperfect_robustOpt}.
The objective function is defined as the weighted sum of the mean and standard deviation of the buckling load so that robust optimisation yields an optimised design that is less sensitive to geometric uncertainties. In robust stability optimisation, the objective is to maximise the mean and minimise the standard deviation of the buckling load, and the balance between the two is controlled by a trade-off parameter. The assumed imperfections are parameterised by a multivariate random vector with a given probability distribution. Consequently, the buckling load of the imperfect structure is also random and has a corresponding probability distribution, which is a push-forward of the probability density of the imperfections~\cite{YelunCirak2024}. As the buckling load is obtained by nonlinear finite element analysis, the probability distribution of the buckling load has no closed-form expression. To estimate the mean and standard deviation of the buckling load, a straightforward approach is to perform Monte Carlo sampling of the given imperfection distribution and to compute the buckling load corresponding to each sample. The two primary techniques for analysing nonlinear structures are the perturbation and path-following methods.  The perturbation method, pioneered by Koiter~\cite{Koiter1967, Budiansky1974}, relies on the asymptotic expansion of the displacement and a load parameter at a known equilibrium configuration.  Its reliance on linear pre-buckling behaviour and small geometric imperfections hinders its application to a broader range of buckling problems. Continuous efforts have been made to improve the perturbation method~\cite{Flores1992, Schafer2006, Liang2013}, which address some of its limitations and allow for broader applications, especially in structural optimisation where repetitive buckling computations are required~\cite{Liguori2018, Vescovini2019}.  In contrast, path-following methods, such as the arc-length method, are more versatile by iteratively computing points along the equilibrium path~\cite{Riks1984}. This iterative process fits seamlessly into the finite element analysis framework and is standard for buckling analysis using finite element analysis packages. However, since stability points are determined only after traversing the equilibrium path, path-following methods can be computationally expensive, making them less suitable for scenarios where many buckling loads must be computed with Monte Carlo sampling. To address this shortcoming, we employ the extended system method~\cite{Wriggers1988}, which computes stability points and the corresponding buckling loads directly without the need to determine the entire equilibrium path. The integration of the stability criterion in the path-following computation allows it to directly converge to stability points, facilitating the analysis of buckling and imperfection sensitivity~\cite{Reitinger1995}.

As mentioned, the most straightforward approach to determine the probability distribution of the buckling load is by Monte Carlo sampling~\cite{Liguori2018, Schevenels2011} from the imperfection probability distribution. This approach is usually computationally expensive since the number of samples can be very large and each sample requires a finite element evaluation. Many techniques can be applied in conjunction with Monte Carlo sampling to improve its efficiency, including Karhunen-Loeve (KL) expansion~\cite{Zhao2014, Chen2010}, expansion optimal linear estimation (EOLE)~\cite{Bai2021, Jansen2013} and perturbation method~\cite{YelunCirak2024, Doltsinis2005, Asadpoure2011}. These methods use truncated series expansions to reduce the dimensionality of the random variable or the number of finite element evaluations; therefore, they have compromised accuracies. We use the quasi-Monte Carlo sampling, i.e. Sobol sampling~\cite{Sobol1967}, to improve sampling efficiency. In contrast to a (pseudo) random sequence, in a quasi-random sequence the points are correlated and are more uniformly distributed, increasing the convergence rate of Monte Carlo sampling from~$\set{O}(N^{-1/2})$ to~$\set{O}( (\log N)^c N^{-1})$, where~$N$ is the number of samples and~$c$ is a constant~\cite{Caflisch_1998,Lemieux2009}. After obtaining the quasi-random sequence, a desired distribution, like the Gaussian distribution, can be obtained by a transformation of the samples.

Conventionally, in structural optimisation, gradient-based methods like the sequential quadratic programming (SQP)~\cite{Schittkowski1986} or the method of moving asymptotes (MMA)~\cite{Svanberg1987} are preferred, and the required sensitivities are determined using an adjoint approach. The same methodology can be extended to the optimisation of buckling loads of structures with geometric nonlinearity~\cite{mroz1994design, zhang2023finite}. Furthermore, perturbation techniques, such as the Koiter method~\cite{Koiter1967}, can be used to account for the effect of imperfections on buckling loads, reducing the use of costly path-following~\cite{mroz1994design, ferrari2025optimization}. In robust optimisation, the buckling load is random and, as noted earlier, the objective function consists of the weighted sum of its mean and standard deviation. Although it is conceivable to extend adjoint and perturbation approaches to robust optimisation of buckling loads, we opt for the gradient-free Bayesian optimisation technique~\cite{Shahriari2016} for the sake of simplicity. Different from other gradient-free metaheuristic approaches like evolutionary algorithms~\cite{Kicinger2005, Wang2020, Bessa2018}, simulated annealing~\cite{Bureerat2008, Hasancebi2002} and particle swarm optimisation~\cite{Zhang2022, Tsiptsis2019}, Bayesian optimisation aims to use a minimum number of function evaluations to converge to the global optimum~\cite{garnett_bayesoptbook_2023}. Consequently, the efficiency of the optimisation procedure is considerably improved.

The outline of this paper is as follows. In Section~\ref{sec:robustOptFormulation} the problem formuation of robust optimisation is defined, and Section~\ref{sec:geoNonlinear} reviews the geometrically nonlinear truss analysis, including the equilibrium equations and tangent stiffness matrix. We present the extended system method in Section~\ref{sec:directRandBasic} for directly determining stability points. The quasi-Monte Carlo sampling is introduced in Section~\ref{sec:directRandMCsampling}, and the algorithm for computing the statistics of buckling loads is proposed in Section~\ref{sec:directRandBucklingDist}. Section~\ref{sec:bayesianOpt} briefly explains the Gaussian process and Bayesian optimisation. In Section~\ref{sec:examples}, we present three examples with increasing complexity to demonstrate the proposed robust optimisation framework.

%% file: geoNonlinear.tex
\section{Problem formulation \label{sec:}}

\subsection{Robust optimisation \label{sec:robustOptFormulation}}
We consider a pin-jointed truss structure in~$\mathbb R^3$ with~$n_p$ nodes, $n_e$ struts and the external loading~$\lambda \vec f$, where~$\vec f \in \mathbb  R^{n_d}$ is a fixed vector, with $n_d$ the number of degrees-of-freedom of the truss, and~$\lambda \in \mathbb R $ is a scaling, or load, parameter. We assume that the truss is suitably supported such that $n_d < 3n_p$ and kinematically stable with no zero energy modes. The equilibrium of the truss is given by
\begin{equation} \label{equ:equilibriumEqu}
	\vec{r}(\vec{x}, \lambda) =  \vec{t}(\vec{x}) - \lambda\vec{f} = \vec{0} \, ,
\end{equation}
where $\vec{r} \in \mathbb R^{n_d}$ is the residual force, $\vec{t} \in \mathbb R^{n_d}$ the internal nodal force, $\vec{x} \in \mathbb R^{n_d}$ the current nodal position vector. The stability of the truss structure is compromised when the external force exceeds a critical buckling load, i.e., a limit or bifurcation load. In stability analysis, the (dimensionless) load parameter $\lambda$ is used to represent the magnitude of the external load. When the structure buckles, the buckling load is denoted by the load parameter $\lambda_c$.

The first few buckling modes typically give the critical imperfection shapes that yield the largest reduction in the buckling~\cite{Budiansky1974, Reitinger1995, Ning2015}. Therefore, in robust stability optimisation with geometric imperfections, the undeformed, or initial, geometry of the truss $\vec{X} \in \mathbb R^{3n_p}$  is considered random and parameterised in terms of the as-designed structure's buckling modes obtained via eigenvalue buckling analysis, i.e., 
\begin{equation} \label{eq:initialConfig}
	\vec{X} = \vec{X}_0 + \sum_{i=1}^{n_b}\beta_i\vec{\phi}_i \, ,
\end{equation}
where $\vec{X}_0 \in \mathbb R^{3n_p}$ is the as-designed geometry, $\vec{\phi}_i \in \mathbb R^{3n_p}$ the $i$-th buckling mode and $n_b$ the number of considered buckling modes. The parameter $\beta_i \in \mathbb R$ is a Gaussian random variable with the probability distribution 
\begin{equation}
	\beta_i \sim \mathcal{N}(\overline{\beta}_i, \sigma_{\beta_i}^2) = \frac{1}{\sqrt{2\pi\sigma_{\beta_i}^2}}\exp\left(-\frac{\left(\beta_i - \overline{\beta}_i\right)^2}{2\sigma_{\beta_i}^2}\right) \, ,
\end{equation}
where $\overline{\beta}_i$  and $\sigma_{\beta_i}$ are the prescribed mean and standard deviation. Note that the parameters~$\beta_i \, (i = 1, \cdots, n_b)$ are assumed uncorrelated.

As the undeformed geometry is random, the corresponding buckling load $\lambda_c$ is also random. The objective of robust optimisation is to maximise the mean $\overline{\lambda}_c$ of the buckling load and simultaneously to minimise its standard deviation~$\sigma_{c}$.  The cross-sectional areas $\vec{a} \in \mathbb R^{n_e}$ of the struts are the design variables. The optimisation is subject to several constraints, specifically, the volume is prescribed with $V_0$, and the cross-sectional areas must lie in the prescribed range $[\vec{a}_{\min}, \, \vec{a}_{\max}]$. The struts carry only axial loads and do not buckle locally, and each strut has a uniform cross-sectional area along its length.

The robust optimisation problem is stated as follows,
\begin{subequations} \label{eq:problem}
	\begin{align}
		\mathop{\text{maximise}}_{\vec{a}} \quad & g(\vec{a}) = \alpha \frac{\overline{\lambda}_c(\vec{a})}{\overline{\lambda}_c^*} - (1-\alpha) \frac{\sigma_{c}(\vec{a})}{\sigma_{c}^*} \, ,\\
		\text{subject to} \quad & 	\vec{r}(\vec{x}, \lambda) =\vec{0} \, , \\
		& V(\vec{a}) = \vec{a} \cdot \vec{l} = V_0 \label{eq:cons2} \, , \\
		& \vec{a}_{\min} \leq \vec{a} \leq \vec{a}_{\max} \, , 
	\end{align}
\end{subequations}
where $\alpha \in [0, 1]$ is a chosen coefficient for a trade-off between the two objectives, and $\overline{\lambda}_c^*$ and $\sigma_{c}^*$ are two normalising constants, i.e. the maximum mean and maximum standard deviation of buckling loads. The two normalising constants~ $\overline{\lambda}_c^*$ and $\sigma_{c}^*$  are determined by performing two different optimisations that consider only maximising the mean $\overline{\lambda}_c$ and standard deviation $\sigma_c$ of buckling loads, respectively. Note that the robust optimisation problem can also be interpreted as a multicriteria optimisation so that the obtained solution for a given~$\alpha$ represents a Pareto point.

\subsection{Geometrically nonlinear analysis \label{sec:geoNonlinear}}

Trusses exhibit large deformations when buckling occurs. Hence, geometric nonlinearities must be considered when determining buckling loads. The formulations of internal force and tangent stiffness matrix for geometrically nonlinear trusses are reviewed briefly in the following~\cite{Bonet2008}.

For a single strut, the Cauchy, or true, stress $\sigma$ is defined in terms of the logarithmic strain $\varepsilon$ as
\begin{equation} \label{eq:cauchyStress}
	\sigma = \frac{V}{v}E\varepsilon = \frac{V}{v}E \ln \left(\frac{l}{L}\right) \, ,
\end{equation}
where $V$ and $v$ are the initial and the deformed volumes, $L$ and $l$ are the initial and the deformed lengths, and $E$ is the Young's modulus.
The internal axial forces at the two end nodes of the strut with the label $e$ are given by
\begin{equation} \label{equ:internalFdefinition}
	\vec{t}^e_a = - \sigma a \vec{n} \, , \quad \vec{t}^e_b = + \sigma a \vec{n} \, .
\end{equation}
Here, $a$ is the deformed cross-sectional area, and $\vec{n}$ is the unit vector representing the direction of the deformed strut. By introducing~\eqref{eq:cauchyStress}, the internal nodal forces can be expressed as
\begin{equation} \label{equ:internalF}
	\vec{t}^e_b = -\vec{t}^e_a = \frac{VE}{l} \ln \left(\frac{l}{L}\right) \vec{n} \, .
\end{equation}
The equilibrium equation of the truss~\eqref{equ:equilibriumEqu} is obtained by equilibrating the internal and external forces at all nodes of the truss, 

 A  path-following approach to compute the nonlinear force-displacement curve adds a constraint condition $h(\vec{x}, \lambda) = 0$ (representing, e.g., a force, a displacement or an arc-length control) to the equilibrium equations. The resulting nonlinear system of equations reads  
\begin{equation} \label{equ:pathFollowingSysEqu}
	\tilde{\vec{r}}(\vec{x}, \lambda) = 
	\begin{pmatrix}
		\vec{r}(\vec{x}, \lambda) \\
		h(\vec{x}, \lambda)
	\end{pmatrix} = \vec{0} \, .
\end{equation}
This equation must be linearised to solve it using the Newton-Raphson method. The linearised equation at iteration~$i$ can be stated in matrix form as
\begin{equation}
\begin{bmatrix}
\vec{K} & -\vec{f} \\
\nabla_{\vec{x}} h & \nabla_\lambda h
\end{bmatrix}_i
\begin{pmatrix}
\Delta \vec u \\
\Delta \lambda
\end{pmatrix}_{i+1}
= -
\begin{pmatrix}
\vec{r}(\vec{x}, \lambda) \\
h(\vec{x}, \lambda)
\end{pmatrix}_i \, ,
\end{equation}
where $\vec{K} = \partial\vec{r}/\partial\vec{x}$ is the tangent stiffness matrix, $\Delta\vec{u}_{i+1}$ is the incremental displacement and $\Delta \lambda_{i+1}$ is the incremental load factor at iteration $i$. The stiffness matrix $\vec{K}$ is assembled from element tangent stiffness matrices $\vec{K}^e$ of struts given by
\begin{align}
	\vec{K}^e = 
		\begin{bmatrix}
			\vec{k}^e & -\vec{k}^e \\
			-\vec{k}^e & \vec{k}^e
		\end{bmatrix} \, , \quad \vec{k}^e = \left(\frac{VE}{l^2}-\frac{2T}{l}\right) \vec{n} \otimes \vec{n} + \frac{T}{l}\vec{I} \, ,
\end{align}
where $\vec{I} \in \mathbb R^{3\times3}$ is the identity matrix, $T = \sigma a$ is the internal axial force, and the length $l$ in the current configuration depends on the displacement.

%
\subsection{Direct computation of stability points \label{sec:directRandBasic}}
%
The extended system method is an elegant approach for computing the stability points without tracing the equilibrium path. First proposed by Wriggers et al.~\cite{Wriggers1988}, the extended system method augments the standard finite element equations~\eqref{equ:equilibriumEqu} with additional constraints restricting solutions to stability points.  One apparent constraint for determining stability points is $\det (\vec{K}) = 0$, which must be the case at all stability points. Since it is not straightforward to linearise $\det(\vec{K})$ directly, it is expedient to consider the condition
\begin{equation} \label{equ:additionalCstr}
	\vec{K}\vec{\phi} = \vec{0} \, .
\end{equation}
This condition is equivalent to $\det (\vec{K}) = 0$ as~$\vec \phi$ is an eigenvector corresponding to an eigenvalue~$0$. To avoid the trivial solution \mbox{$\vec{\phi} = \vec{0}$}, an additional constraint on the norm of eigenvector~$\vec \phi$ denoted by $s(\vec{\phi})$ is introduced, that is,
\begin{equation} \label{equ:eigenvectorCstr}
	s(\vec{\phi}) = \|\vec{\phi}\| - 1 = 0 \, .
\end{equation}

Combining equations~\eqref{equ:equilibriumEqu},~\eqref{equ:additionalCstr} and~\eqref{equ:eigenvectorCstr}, the new extended system of equations $\hat{\vec{r}}$ for computing stability points can be written as
\begin{equation} \label{equ:extendedSysOfEqu}
	\hat{\vec{r}}(\vec{x}, \vec{\phi}, \lambda) =  
	\begin{pmatrix}
		\vec{r}(\vec{x}, \lambda) \\
		\vec{K}(\vec{x}, \lambda) \vec{\phi} \\
		s(\vec{\phi})
	\end{pmatrix} = \vec{0} \, .
\end{equation}
Similar to conventional nonlinear finite element analysis, the system of equations~\eqref{equ:extendedSysOfEqu} is linearised and solved with the Newton-Raphson method. The linearisation of~\eqref{equ:extendedSysOfEqu} at iteration $i$ yields
\begin{subequations}
	\begin{align}
		\hat{\vec{K}}_i \Delta \vec{w}_{i+1} &= - \hat{\vec{r}}_i \, , \label{equ:extendedLinearised} \\
		\vec{w}_{i+1} &= \vec{w}_{i} + \Delta \vec{w}_{i+1} \, , 
	\end{align}	
\end{subequations}
where $\hat{\vec{K}}$ is the Jacobian of $\hat{\vec{r}}$, $\vec{w} = (\vec{x}, \, \vec{\phi}, \, \lambda)^\trans$ the vector of unknown variables, and $\Delta\vec{w} = (\Delta\vec{u}, \, \Delta\vec{\phi}, \, \Delta\lambda)^\trans$  the incremental vector. For a truss with $n_d$ degrees-of-freedom, $\vec{w} \in \mathbb R^{2n_d+1} $ contains the nodal positions $\vec{x} \in \mathbb{R}^{n_d}$, eigenvector $\vec{\phi} \in \mathbb{R}^{n_d}$ and load parameter $\lambda \in \mathbb{R}$. When written in matrix form, equation~\eqref{equ:extendedLinearised} reads
\begin{equation} \label{equ:extendedLinearisedDetail}
	\begin{bmatrix}
		\vec{K} & \vec{0} & - \vec{f} \\
		\nabla_{\vec{x}}(\vec{K}\vec{\phi}) & \vec{K} & \nabla_{\lambda}(\vec{K}\vec{\phi}) \\
		\vec{0} & \nabla_{\vec{\phi}}s & 0
	\end{bmatrix}_i
	\begin{pmatrix}
		\Delta \vec{u} \\
		\Delta \vec{\phi} \\
		\Delta \lambda
	\end{pmatrix}_{i+1} = -
	\begin{pmatrix}
		\vec{r}(\vec{x}, \lambda) \\
		\vec{K}(\vec{x}, \lambda) \vec{\phi} \\
		s(\vec{\phi})
	\end{pmatrix}_i \, .
\end{equation}

A partitioning method is used to efficiently solve the system of equations~\eqref{equ:extendedLinearisedDetail}. The incremental displacement $\Delta\vec{u}$ is obtained from the first row in~\eqref{equ:extendedLinearisedDetail} as 
\begin{equation} \label{eq:updateU}
	\Delta \vec{u}_{i+1} = \left(\vec{K}_i^{-1}\vec{f}\right) \Delta \lambda_{i+1} - \vec{K}_i^{-1}\vec{r}_i \, .
\end{equation}
Similarly, the incremental eigenvector $\Delta\vec{\phi}_{i+1}$ is obtained from the second row in~\eqref{equ:extendedLinearisedDetail} as
\begin{equation}\label{equ:extendedDeltaPhi}
	\Delta\vec{\phi}_{i+1} = \Delta\vec{\phi}^{(1)}\Delta\lambda_{i+1} + \Delta\vec{\phi}^{(2)} \, , 
\end{equation}
with
\begin{subequations} \label{equ:extendedPhi12}
\begin{align}
	&\Delta\vec{\phi}^{(1)} = -\vec{K}_i^{-1}\left[\nabla_{\vec{x}}(\vec{K}\vec{\phi})_i \vec{K}_i^{-1}\vec{f} + \nabla_{\lambda}(\vec{K}\vec{\phi})_i\right] \, ,\\
	&\Delta\vec{\phi}^{(2)} = -\vec{K}_i^{-1}\left[(\vec{K}\vec{\phi})_i - \nabla_{\vec{x}}(\vec{K}\vec{\phi})_i \vec{K}_i^{-1}\vec{r}_i\right] \, ,
\end{align}
\end{subequations}
From the third row in~\eqref{equ:extendedLinearisedDetail}, the incremental load parameter $\Delta\lambda_{i+1}$ is  determined as 
\begin{equation} \label{equ:extendedDeltaLambda}
	\Delta\lambda_{i+1} = - \frac{\nabla_{\vec{\phi}}s_i \Delta\vec{\phi}^{(2)} + s_i}{\nabla_{\vec{\phi}}s_i \Delta\vec{\phi}^{(1)}} \, .
\end{equation}
Finally, the unknown variables at iteration $(i+1)$ are given by
\begin{equation}
\vec{x}_{i+1} = \vec{x}_i + \Delta\vec{u}_{i+1} \, , \quad
\vec{\phi}_{i+1} = \vec{\phi}_i + \Delta\vec{\phi}_{i+1} \, , \quad
\lambda_{i+1} = \lambda_i + \Delta\lambda_{i+1} \, .
\end{equation}

\begin{figure*}
	\centering
	\subfloat[Von Mises truss] {
		\hspace{6mm}
		\includegraphics[width=0.37 \textwidth]{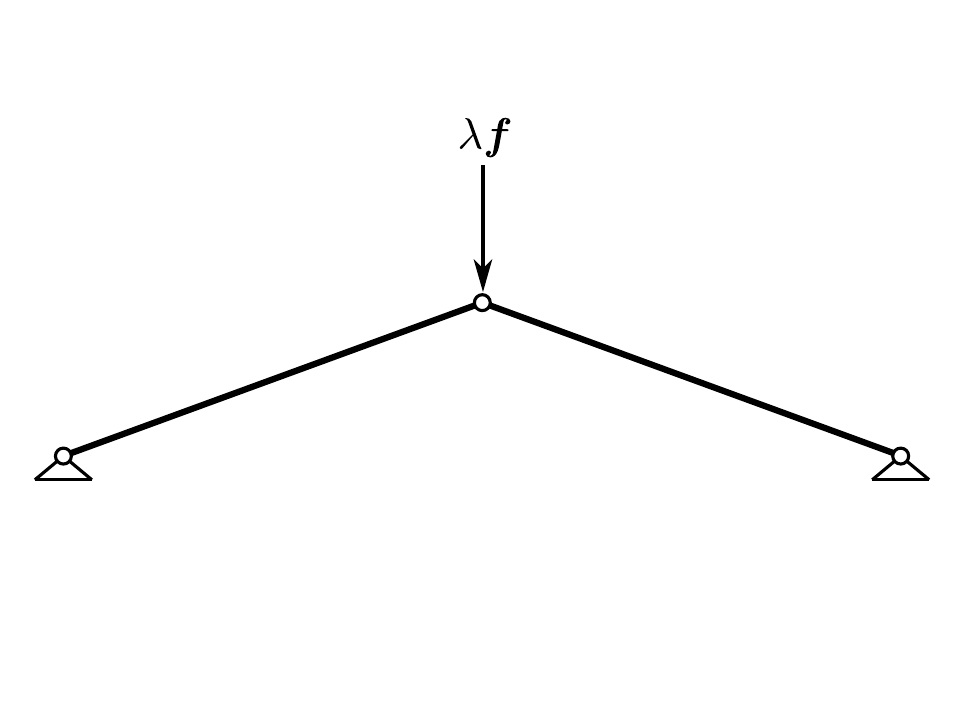}
		\label{fig:vonMises_setup} }
	\hspace{8mm}
	\subfloat[First stability point] {
		\centering
		\includegraphics[width=0.4 \textwidth]{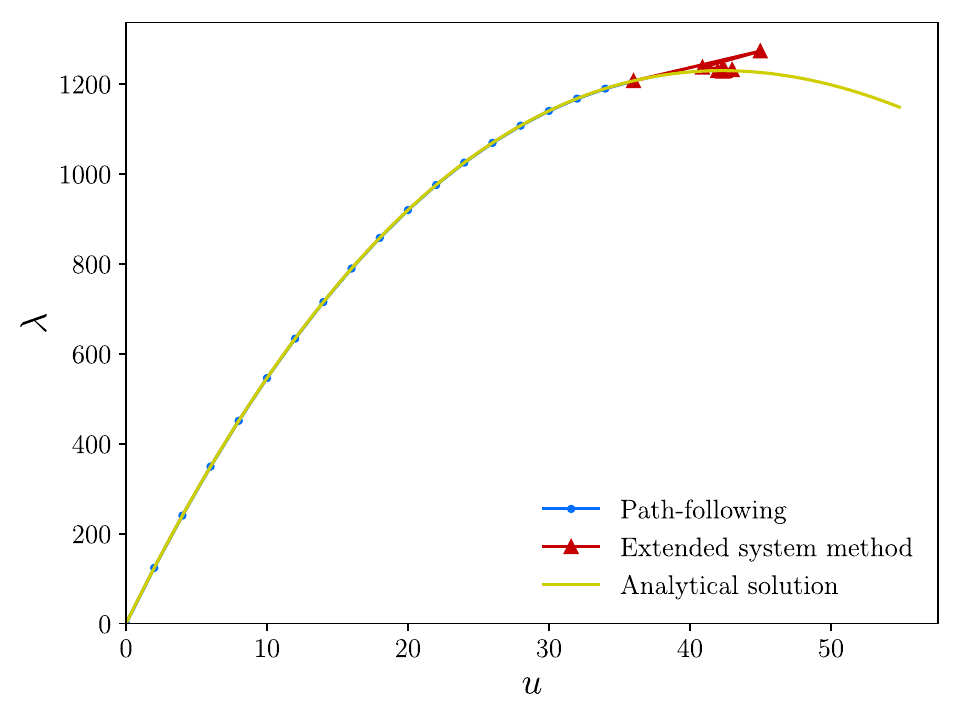}
		\label{fig:vonMisesAnalyticalNLFEAextendedSystem} }
	\\
	\subfloat[First stability point using the extended system method] {
		\includegraphics[width=0.415 \textwidth]{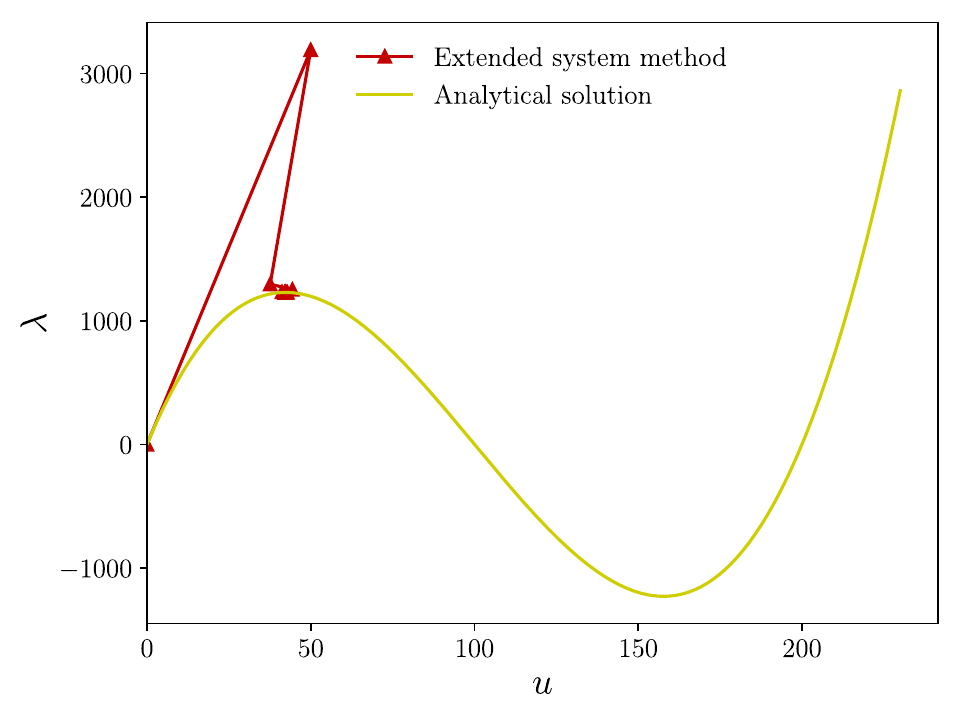}
		\label{fig:vonMises_OriginExtendedSystem_stabilityPt1} }
	\hspace{5mm}
	\subfloat[Second stability point using the extended system method] {
		\includegraphics[width=0.415 \textwidth]{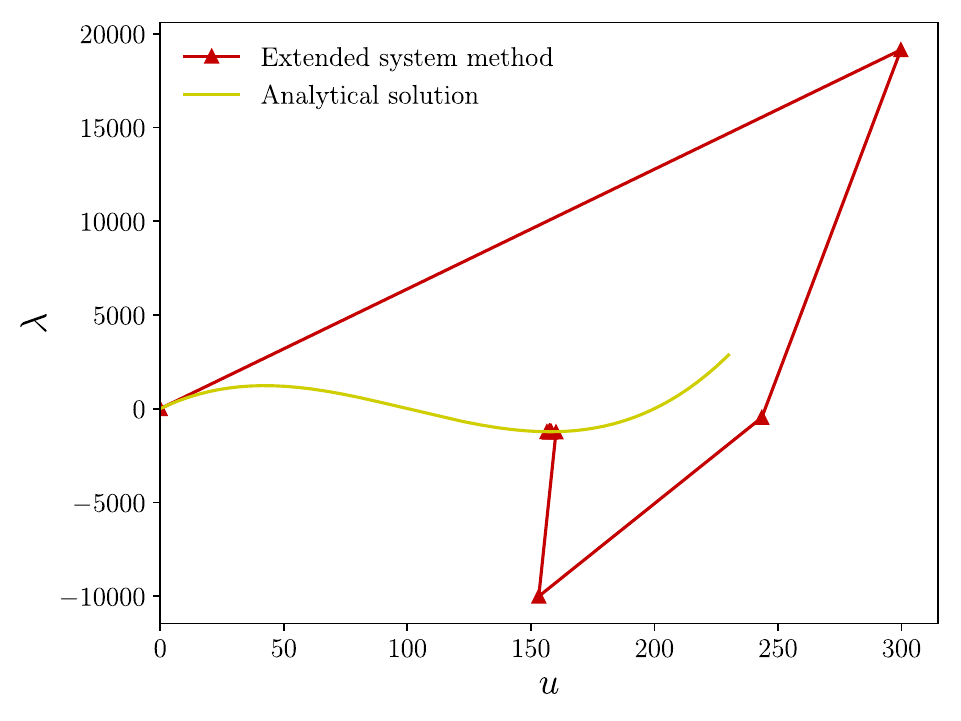}
		\label{fig:vonMises_OriginExtendedSystem_stabilityPt2} }
	\caption{Computation of stability points of a Von Mises truss using the extended system method, with and without path-following. The yellow curves represent the analytical solution of the load-displacement relationship. Figure (b) shows the combination of the nonlinear path-following and extended system method converging to the first stability point. Figures (c) and (d) are the convergence of the extended system method to the first and second stability points, respectively.}
	\label{fig:vonMises_OriginExtendedSystem}
\end{figure*}

The predictor values of the three variables are required at the beginning of iteration~$i=0$. Correctly choosing the initial eigenvector $\vec{\phi}_0$ can lead to a fast convergence. In addition, for problems with multiple stability points, to which stability point the iteration converges depends on the choice of $\vec{\phi}_0$, see Figure~\ref{fig:vonMises_OriginExtendedSystem_stabilityPt1} and~\ref{fig:vonMises_OriginExtendedSystem_stabilityPt2}. The choice of $\vec{\phi} $ is usually inconsequential when first an equilibrium point close to the stability point is found using path-following, and then the extended system is considered, see Figure~\ref{fig:vonMisesAnalyticalNLFEAextendedSystem}.

One of the challenges in the extended system method is the computation of the directional derivatives of the tangent stiffness matrix appearing in~\eqref{equ:extendedPhi12}. Note that for conservative external loading, the term $\nabla_{\lambda}(\vec{K}\vec{\phi})$ is zero and the directional derivatives can be derived analytically for truss structures. Nevertheless, for more complicated structures, e.g. shells, the derivatives may need to be approximated numerically. Apart from the finite difference approximation in the original papers by Wriggers et al.~\cite{Wriggers1988, Wriggers1990}, the complex-step derivative method~\cite{Lyness1967, Squire1998, Martins2003, Lai2008} provides more accurate and reliable results regardless of the choice of the perturbation parameter.

%% file: directRandStability.tex
\section{Statistics of random buckling loads \label{sec:directRandStability}}
The mean and standard deviation of buckling loads are estimated by sampling random imperfections and computing their corresponding buckling loads. A quasi-Monte Carlo approach is used to generate the random imperfections, and the extended system method is employed to compute the buckling load for each configuration.

\subsection{Quasi-Monte Carlo sampling  \label{sec:directRandMCsampling}} 

In contrast to Monte Carlo sampling, in quasi-Monte Carlo sampling a low-discrepancy quasi-random sequence is employed~\cite{Robert2004}. Since the quasi-random sequence has more uniformly distributed points than a random sequence, the quasi-Monte Carlo sampling has a faster convergence rate of $\mathcal{O}((\log N)^c N^{-1})$~\cite{Caflisch_1998, Lemieux2009}, compared to $\mathcal{O}(N^{-1/2})$ for the conventional Monte Carlo sampling, and gives better accuracy with a smaller number of samples.

\begin{figure}
	\centering
	\subfloat[$1-128$ points] {
		\includegraphics[width=0.3 \textwidth]{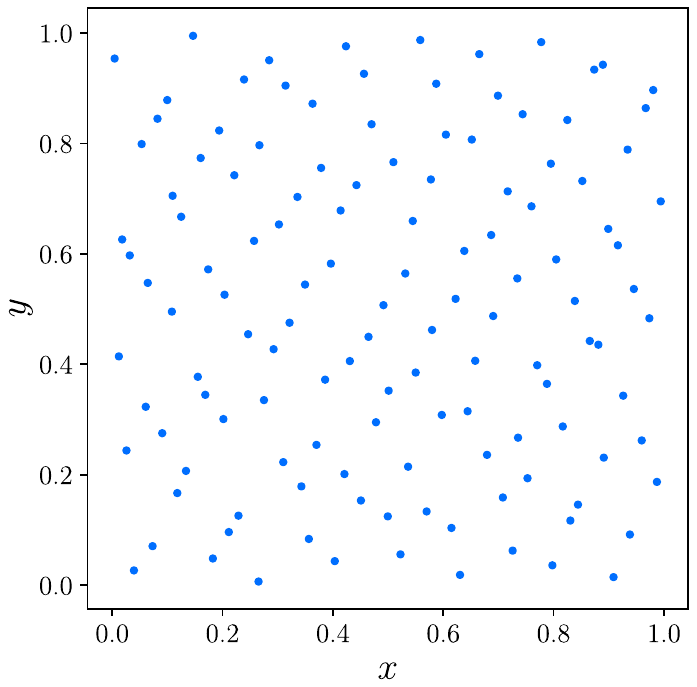}
		\label{fig:sobolSequence2dset1} }
	\hfill
	\subfloat[$129-256$ points] {
		\includegraphics[width=0.3 \textwidth]{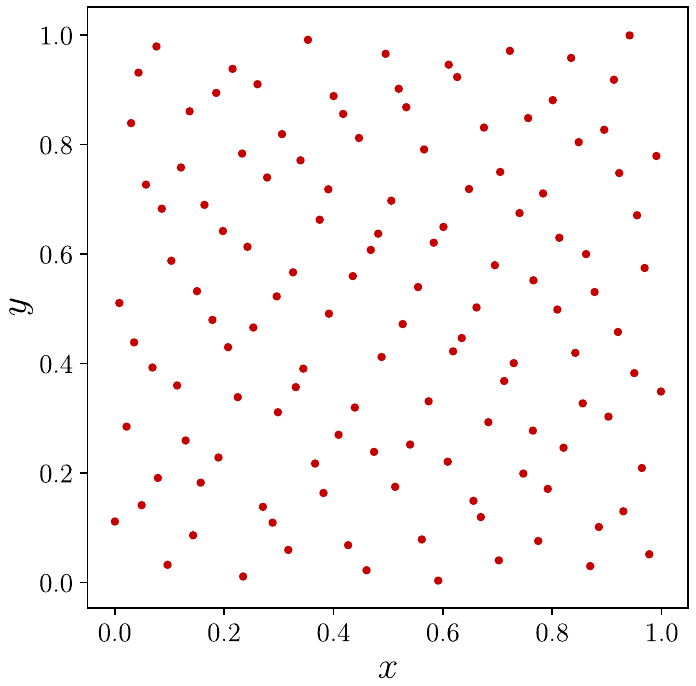}
		\label{fig:sobolSequence2dset2} }
	\hfill
	\subfloat[$1-256$ points] {
		\includegraphics[width=0.3 \textwidth]{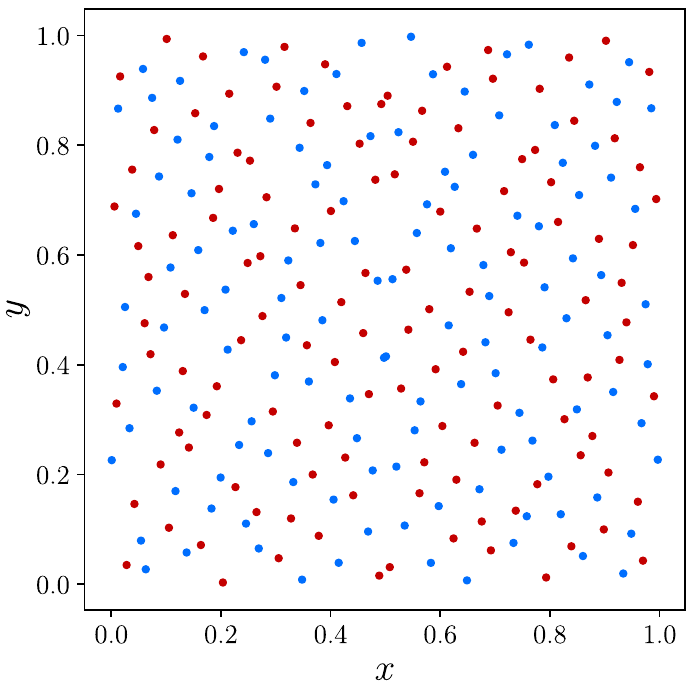}
		\label{fig:sobolSequence2dsetAll} }
	\caption{A $2$D Sobol sequence with consecutive sets. Figures (a) and (b) show two consecutive sets with $2^7=128$ points in each set. Figure (c) shows the combination of the two sets with $256$ points. Note that each set is fairly uniformly distributed.}
	\label{fig:sobolSequence2d}
\end{figure}

In our computations, the Sobol sampling using a quasi-random sequence, also called a Sobol sequence~\cite{Sobol1967}, is employed to sample from the probability distribution of geometric imperfections. It generates a set of sample points that are more evenly distributed than Monte Carlo sampling. The Sobol sequence is suitable in large dimensions~\cite{Lemieux2009} and its generation mainly requires bitwise operations~\cite{Bratley1988, Joe2003, Antonov1979}. With suitable choices of direction numbers, every set with $2^s$ points of Sobol sequence is well distributed in the $s$-dimensional hypercube and the combination of successive sets also satisfies the low-discrepancy property. Figure~\ref{fig:sobolSequence2d} depicts a 2D Sobol sequence as an example, in which there is no overlap between the two sets of points in Figures~\ref{fig:sobolSequence2dset1} and~\ref{fig:sobolSequence2dset2}. The space is successively refined as the Sobol sequence progresses, guaranteeing a more even space filling. Hence, the Sobol sequence can be used to generate samples in a set-wise manner. In this work, the Python library \texttt{SciPy} is used to generate the Sobol sequence.

The geometric imperfection in practice has a different probability density other than a uniform density, e.g. a Gaussian density. The uniform density of the Sobol sequence has to be transformed into other specified probability densities. Inverse transform sampling using the CDF (cumulative distribution function) is adopted to transform uniformly distributed samples to samples from the desired probability density. To illustrate the effectiveness of quasi-random sequence in reducing the number of samples, the comparison of the normally-distributed random and quasi-random Sobol samples is shown in Figure~\ref{fig:vonMises_imperfPct_imperfBucklingF_differentDist}. Having the same number of samples, the quasi-random samples in Figure~\ref{fig:vonMises_imperfPct_normal} fit the normal distribution far better than the random samples in Figure~\ref{fig:vonMises_imperfPct_normal_rand}; in other words, the number of samples can be significantly reduced with quasi-random sampling.

\subsection{Buckling load statistics \label{sec:directRandBucklingDist}}
The random initial configuration $\vec{X}$, cf.~\eqref{eq:initialConfig}, is obtained via quasi-Monte Carlo sampling of the imperfections. First, a Sobol sequence $\vec \gamma = \{\vec\gamma^{(i)}\}_{i=1, \cdots, 2m}$ following the uniform density is computed, where $\vec\gamma^{(i)} = (\gamma^{(i)}_1, \cdots, \gamma^{(i)}_{n_b})^\trans$. Using the inverse transform sampling, it is then transformed to the amplification factors $\vec \beta = \{\vec\beta^{(i)}\}_{i=1, \cdots, 2m}$ of buckling modes $\vec\Phi = (\vec\phi_1, \cdots, \vec\phi_{n_b})^\trans$, where $\vec\beta^{(i)} = (\beta^{(i)}_1, \cdots, \beta^{(i)}_{n_b})^\trans$ with each component following the normal distribution that is assumed to have a zero mean and a fixed standard deviation, i.e., $\beta^{(i)}_k \sim \mathcal{N}(0, \sigma_\beta^2)$. The sampled amplification factors $\vec \beta$ are separated into positive and negative ones denoted as $\vec{\beta}^{+}$ and $\vec{\beta}^{-}$, respectively. The positive imperfection factors are sorted in ascending order and the negative ones in descending order. The buckling loads $\lambda_c$ are computed for $\vec{\beta}^{+}$ and $\vec{\beta}^{-}$ separately.

\begin{figure}
	\centering
	\subfloat[Positive imperfect structures] {
		\includegraphics[width=0.45 \textwidth]{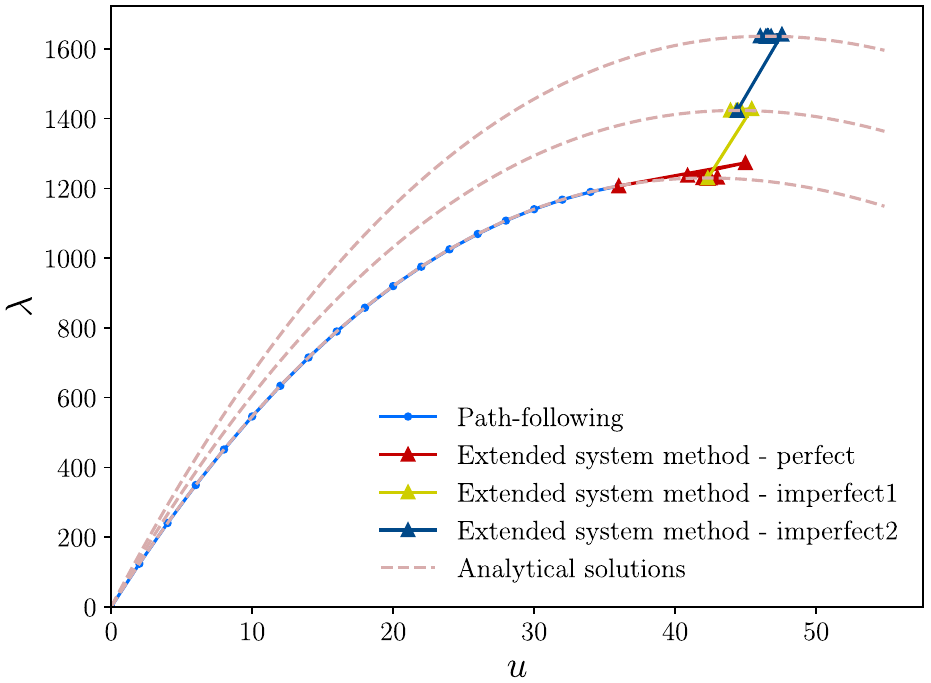}
		\label{fig:vonMises_imperfBucklingF_algorithm_posi} }
	\hspace{5mm}
	\subfloat[Negative imperfect structures] {
		\includegraphics[width=0.45 \textwidth]{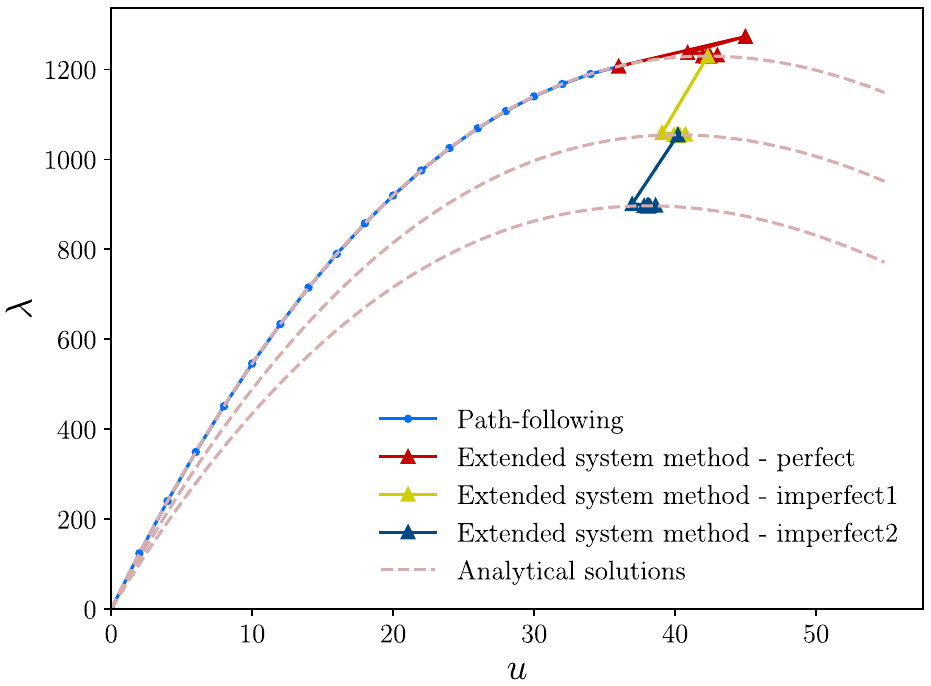}
		\label{fig:vonMises_imperfBucklingF_algorithm_nega} }
	\caption{Illustration of determining the buckling loads of a von Mises truss with random geometric imperfections. Positive and negative amplification factors are separated and sorted. The repeated computations start with the as-designed truss using the combination of the path-following and the extended system method and then iterate through every geometric imperfection with the extended system method only.}
	\label{fig:vonMises_imperfBucklingF_algorithm}
\end{figure}

For either $\vec{\beta}^+$ or $\vec{\beta}^-$, we start from the as-designed configuration $\vec{X}_0$ and compute its buckling load $\lambda_c^0$. To ensure the convergence to its first stability point, the path-following method is adopted at the beginning of iterations and it is switched to the extended system method when $\det (\vec{K})$ becomes sufficiently small or negative. Starting from $\lambda_c^0$ we iteratively compute the buckling loads of trusses with random configurations $\vec X$ using the extended system method only. The stability point of a configuration serves as the predictor for the iterations for the next configuration, as illustrated in Figure~\ref{fig:vonMises_imperfBucklingF_algorithm}. The process of computing the empirical, or sample, mean and standard deviation of the buckling load is summarised in Algorithm~\ref{alg:bucklingFDist}.

\begin{algorithm*} \footnotesize
	\caption{Computation of buckling load statistics  \label{alg:bucklingFDist}}
	
	\begin{algorithmic}[1]
		\State \textbf{Input:} buckling modes $\vec{\Phi} = (\vec\phi_1, \cdots, \vec\phi_{n_b})^\trans$ of the as-designed truss, buckling loads set $\vec{\lambda}_c = \emptyset$ for all the configurations considered, imperfection distribution $\mathcal{N}$, number of samples $2m$
		
		\State $\{\vec \gamma^{(1)}, \dots, \vec \gamma^{(2m)}\} \Leftarrow$ Sample $\vec \gamma \sim\mathcal{U}(0, 1)$ from Sobol sequence
		\State $\{\vec \beta^{(1)}, \dots, \vec \beta^{(2m)}\} \Leftarrow$ Transform $\vec \gamma \sim\mathcal{U}(0, 1)$ to $\vec \beta \sim \mathcal{N}$ with the inverse transform method
		\State $\{(\vec \beta^{+})^{(1)}, \dots, (\vec \beta^{+})^{(m)}\}$ and $\{(\vec \beta^-)^{(1)}, \dots, (\vec \beta^-)^{(m)}\} \Leftarrow$ Separate and sort $\vec \beta$ \Comment{$\vec \beta^+$ in ascending, $\vec \beta^-$ in descending}
		\State Compute geometric imperfections $(\vec{\beta}^+)^{(j)}\cdot\vec{\Phi}$ and $(\vec{\beta}^-)^{(j)}\cdot\vec{\Phi}$ ($j = 1, \cdots, m$)
		\State $\vec{\lambda}_c \gets \vec{\lambda}_c \cup \{\lambda_c^{(0)}\} \Leftarrow$ Analyse the as-designed truss by path-following methods and extended system method
		\For{$i \in \{0, 1\}$} 
			\State Restore to the as-designed configuration $\vec{X}_0$	
			\For{$j \in \{1, \dots m\}$} \Comment{Iterate through imperfect structures; set $(\vec{\beta}^+)^{(0)} = \vec{0}$ and $(\vec{\beta}^-)^{(0)} = \vec{0}$}
				\If{$i == 0$} 
					\State Add $(\vec{\beta}^+)^{(j)}\cdot\vec\Phi - (\vec{\beta}^+)^{(j - 1)}\cdot\vec\Phi$ to the previous configuration
				\ElsIf{$i == 1$}
					\State Add $(\vec{\beta}^-)^{(j)}\cdot\vec\Phi - (\vec{\beta}^-)^{(j - 1)}\cdot\vec\Phi$ to the previous configuration
				\EndIf
				\State $k = j + i\cdot m$
				\State $\lambda_c^{(k)} \Leftarrow$ Analyse imperfect structure by extended system method
				\State $\vec\lambda_c \gets \vec\lambda_c \cup \{\lambda_c^{(k)}\}$
			\EndFor
		\EndFor
		\State \textbf{return:} $\vec\lambda_c$
	\end{algorithmic}
	
\end{algorithm*}

\begin{figure}
	\centering
	\subfloat[Imperfections with normal distribution (quasi-random)] {
		\includegraphics[width=0.45 \textwidth]{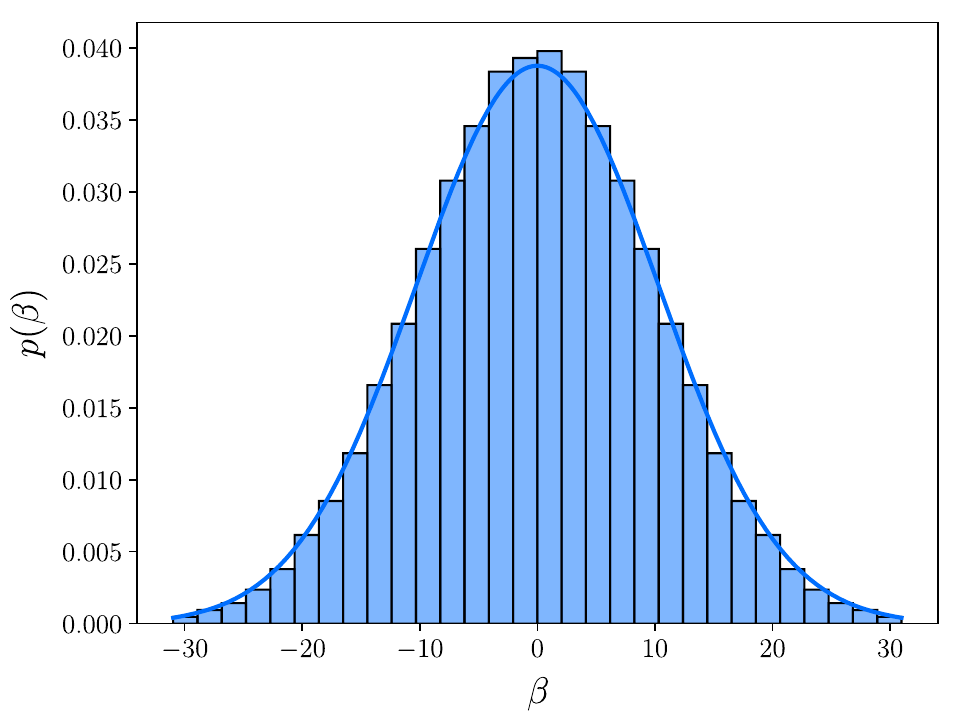}
		\label{fig:vonMises_imperfPct_normal} }
	\hfil	
	\subfloat[Buckling loads (quasi-random)] {
		\includegraphics[width=0.45 \textwidth]{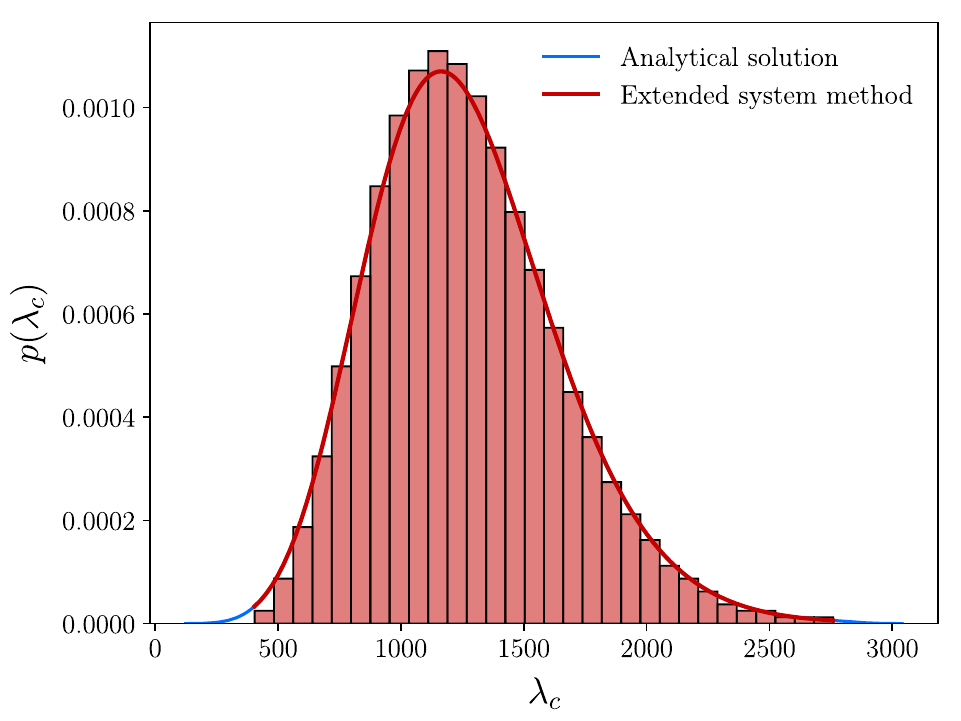}
		\label{fig:vonMises_imperfBucklingF_normal} }
	\\
	\subfloat[Imperfections with normal distribution (pseudorandom)] {
		\includegraphics[width=0.45 \textwidth]{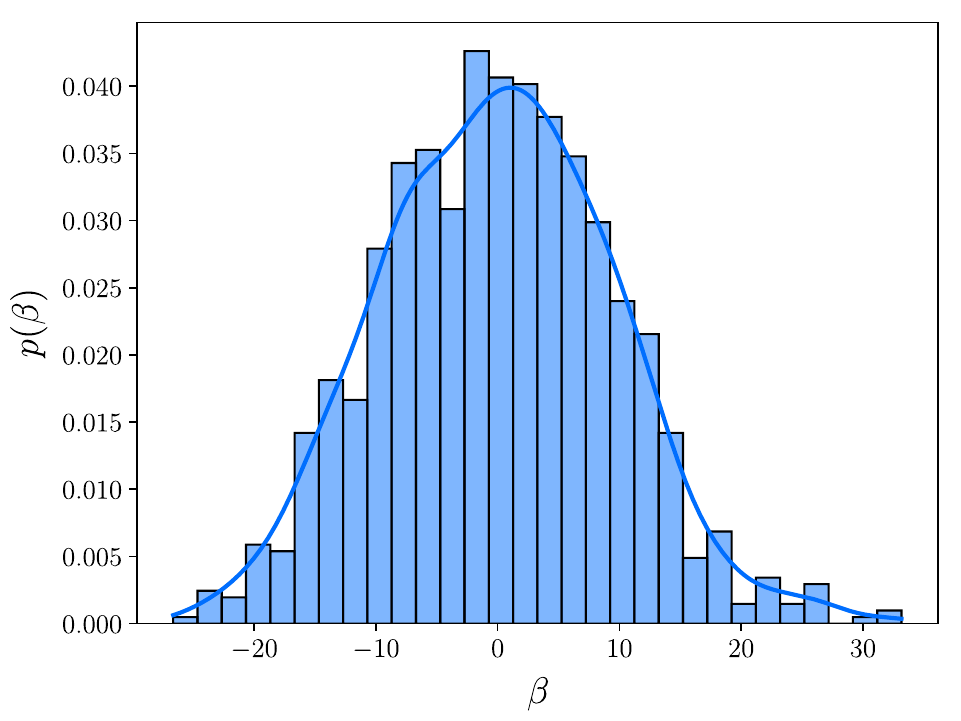}
		\label{fig:vonMises_imperfPct_normal_rand} }
	\hfil	
	\subfloat[Buckling loads (pseudorandom)] {
		\includegraphics[width=0.45 \textwidth]{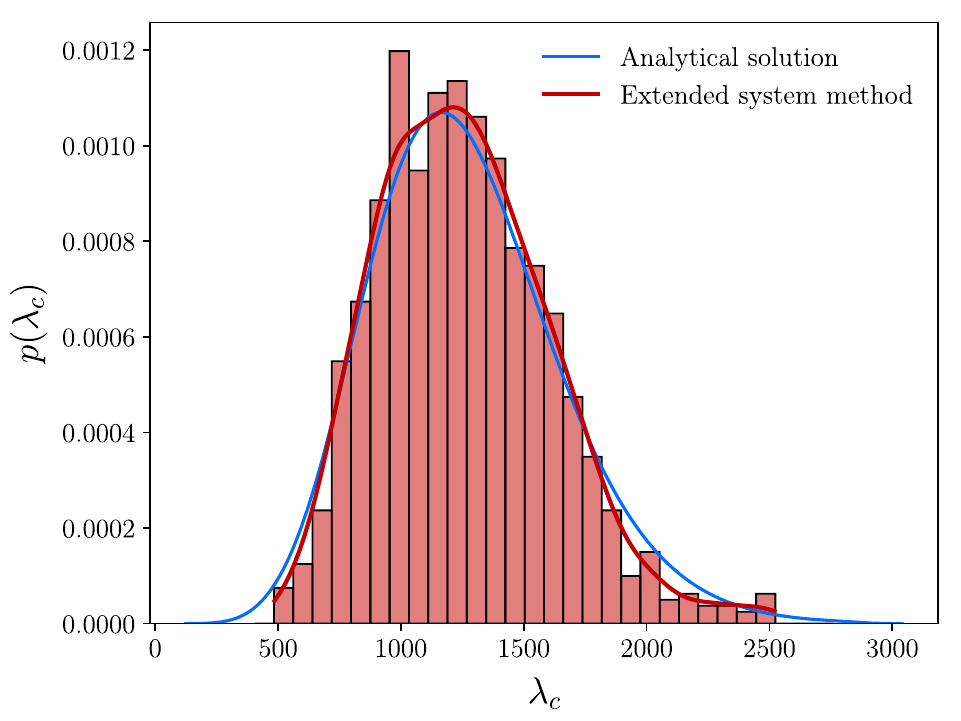}
		\label{fig:vonMises_imperfBucklingF_normal_rand} }
	\caption{Comparison of the quasi-random Sobol sampling and pseudorandom sampling of the geometric imperfection parameter $\beta$ and the corresponding histograms of the buckling load $\lambda_c$. All histograms use 1024 samples. The solid lines are kernel density estimates.}
	\label{fig:vonMises_imperfPct_imperfBucklingF_differentDist}
\end{figure}

The buckling load density computed with the proposed algorithm is compared with that by sampling from the analytical solution, using a von Mises truss with a normal density for the random geometric imperfection. A total of $1024$ quasi-random samples of $\vec \beta$ are generated from the Sobol sequence. With the geometric imperfection distribution in Figure~\ref{fig:vonMises_imperfPct_normal}, the corresponding buckling load distribution is shown in Figure~\ref{fig:vonMises_imperfBucklingF_normal}. The buckling load histograms obtained from the proposed algorithm and analytical solutions are almost identical. Compared with the buckling load histogram (Figure~\ref{fig:vonMises_imperfBucklingF_normal_rand}) using pesudorandom samples in Figure~\ref{fig:vonMises_imperfPct_normal_rand}, the one using quasi-random samples is much closer to the analytical solution. Therefore, in addition to the number of samples, the number of finite element evaluations can also be significantly reduced with quasi-random sampling.

\begin{figure}
	\centering
	\subfloat[Empirical mean of buckling load] {
		\includegraphics[width=0.45 \textwidth]{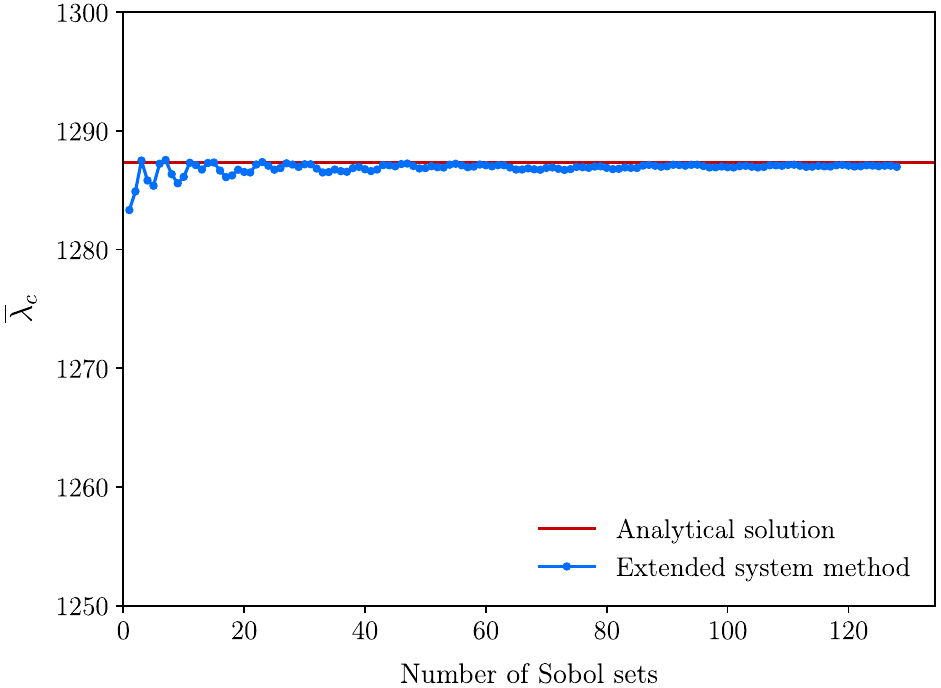}
		\label{fig:vonMises_analyticalNumerical_convergence_mean} }
	\hspace{5mm}
	\subfloat[Standard deviation of buckling load] {
		\includegraphics[width=0.45 \textwidth]{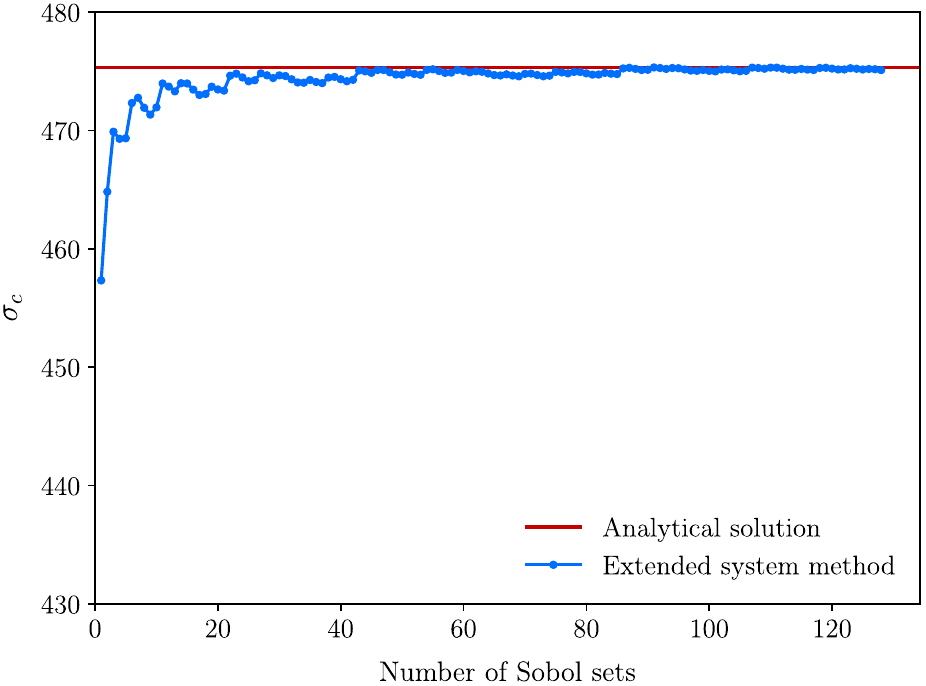}
		\label{fig:vonMises_analyticalNumerical_convergence_std} }
	\caption{Convergence of the empirical mean and standard deviation of the buckling load of a von Mises truss, using an increasing number of sample sets.}
	\label{fig:vonMises_analyticalNumerical_convergence}
\end{figure}

The convergence of the algorithm is demonstrated with the empirical mean and standard deviation of the buckling load obtained with an increasing number of samples, as shown in Figure~\ref{fig:vonMises_analyticalNumerical_convergence}. The geometric imperfections follow a normal distribution and each sample set contains $128$ quasi-random samples. The analytical solution is obtained using $16384$ random samples. It can be seen that the mean and standard deviation values converge to the analytical ones quickly within a small number of sample sets. Furthermore, even with one set of $128$ samples, the errors between our results and the analytical result are very small, $0.3\%$ in mean and $3.7\%$ in standard deviation. Therefore, it is sufficient to use a small number of quasi-random samples (e.g. $128$ samples) in our algorithm to compute the mean and standard deviation of the buckling load with reasonable accuracy, resulting in a significant reduction in computational cost.

%% file: bayesianOpt.tex
\section{Surrogate-based robust optimisation} \label{sec:bayesianOpt}
%
We solve the robust optimisation problem~\eqref{eq:problem} using Bayesian optimisation. The Gaussian process (GP) model serves as the surrogate model of the objective function~\cite{Archetti2019}. The optimised cross-sectional areas of struts are obtained by iteratively adding new data points that maximise an acquisition function to the training dataset of the GP model. To enforce the volume constraint~\eqref{eq:cons2}, one of the cross-sectional areas, for example the last design variable $a_{n_c}$, is expressed in terms of the other cross-sectional areas in the objective function, such that the linear equality constraint is embedded in optimisation and the variable $a_{n_c}$ is obtained straightforwardly after other design variables are optimised.

\subsection{Gaussian process surrogate}

A Gaussian process represents a collection of random variables, and any finite number of those random variables have a joint Gaussian distribution~\cite{Rasmussen2004}. A Gaussian process is defined by a mean function $\overline{f}(\vec{a})$ and a covariance function $\mathrm{cov}(\vec{a}, \vec{a}')$ as $f(\vec{a}) \sim \mathcal{GP}(\overline{f}(\vec{a}), \mathrm{cov}(\vec{a}, \vec{a}'))$, where
\begin{subequations}
	\begin{align}
		& \overline{f}(\vec{a}) = \expect\left[f(\vec{a})\right] \, , \\
		& \mathrm{cov}(\vec{a}, \vec{a}') = \expect\left[(f(\vec{a})-\overline{f}(\vec{a}))(f(\vec{a}')-\overline{f}(\vec{a}'))\right] \, ,
	\end{align}
\end{subequations}
with $\expect$ denoting the expectation operator. We choose the Mat\'ern covariance function given by
\begin{equation} \label{eq:maternCovariance}
	\mathrm{cov}(\vec{a}, \vec{a}') = \frac{2^{1-\nu}}{\Gamma(\nu)} \left(\frac{\sqrt{2\nu} }{\eta}r\right)^{\nu} K_{\nu} \left(\frac{\sqrt{2\nu} }{\eta}r\right) \, ,
\end{equation}
where $r=\|\vec{a}-\vec{a}'\|$, $\nu > 0$ is the smoothness parameter, $\eta > 0$ is the lengthscale parameter, $\Gamma$ is the Gamma function and $K_{\nu}$ is the modified Bessel function.

After fitting with a given training dataset $\{\vec{A}, \vec{g}\} = \{(\vec{a}_i, g(\vec{a}_i)) \}_{i=1}^n$, the Gaussian process regression can predict the mean and variance at the given $n_*$ test points $\vec{A}_*$. The statistical observation model assumed in Gaussian process regression reads 
\begin{equation}
	g(\vec a) = f(\vec a) + \epsilon \, ,
\end{equation}
where~$\epsilon \sim \mathcal{N}(0, \sigma_\epsilon^2)$ is the observation noise, and~$f(\vec a) $ and~$\epsilon$ are statistically independent.  The joint probability distribution of the training outputs $\vec{g}$ and test outputs $\vec{g}_*$ is given by
\begin{equation}
	\begin{bmatrix}
		\vec g \\
		\vec g_*
	\end{bmatrix}
	\sim \mathcal{N}\left(
	\vec{0}  ,
	\begin{bmatrix}
		\vec{C}(\vec{A}, \vec{A}) + \sigma_\epsilon^2 \vec{I} & \vec{C}(\vec{A}, \vec{A}_*) \\
		\vec{C}(\vec{A}_*, \vec{A}) & \vec{C}(\vec{A}_*, \vec{A}_*)
	\end{bmatrix}\right) \, ,
\end{equation}
where $\vec{C}(\cdot, \cdot) \in \mathbb{R}^{n \times n}$ is the covariance matrix with each entry evaluated by plugging the training and test points in the covariance function~\eqref{eq:maternCovariance}, and $\vec{I}$ is the identity matrix. The predictive distribution of test outputs $\vec{g}_*$ for given $\vec{A}_*$ and the training dataset $\{\vec{A}, \vec{g}\}$ is given by
\begin{equation}
	\vec{g}_* \vert \vec{A}_*, \vec{A}, \vec{g} \sim \mathcal{N}\left(\overline{\vec{g}}_*, \mathrm{cov}(\vec{g}_*)\right) \, ,
\end{equation}
with the mean vector and covariance matrix 
\begin{subequations}
	\begin{align}
		&\overline{\vec{g}}_* = \vec{C}(\vec{A}_*, \vec{A}) \left[\vec{C}(\vec{A}, \vec{A}) + \sigma_\epsilon^2 \vec{I}\right]^{-1} \vec{g} \, , \\
		&\mathrm{cov}(\vec{g}_*) = \vec{C}(\vec{A}_*, \vec{A}_*)  - \vec{C}(\vec{A}_*, \vec{A}) \left[\vec{C}(\vec{A}, \vec{A}) + \sigma_\epsilon^2 \vec{I}\right]^{-1} \vec{C}(\vec{A}, \vec{A}_*) \, .
	\end{align}
\end{subequations}

The fitting, or training, of the Gaussian process is achieved by optimising the hyperparameters $\vec{\theta}=(\nu, \, \eta, \, \sigma_{\epsilon})$ through maximising the log marginal likelihood, given by
\begin{equation}
	\log p(\vec{g} | \vec{A}, \vec{\theta}) = -\frac{1}{2} \vec{g}^T \left(\vec{C}+\sigma_{\epsilon}^2 \vec{I} \right)^{-1} \vec{g} - \frac{1}{2} \log \left| \vec{C}+\sigma_{\epsilon}^2 \vec{I} \right| - \frac{n}{2} \log 2\pi \, ,
\end{equation}
where $n$ is the dimension of $\vec{C}$, and $\left| \vec{C}+\sigma_{\epsilon}^2 \vec{I} \right|$ is the determinant.

\subsection{Bayesian optimisation}

\begin{figure}
	\centering
	\includegraphics[height=0.35 \textwidth]{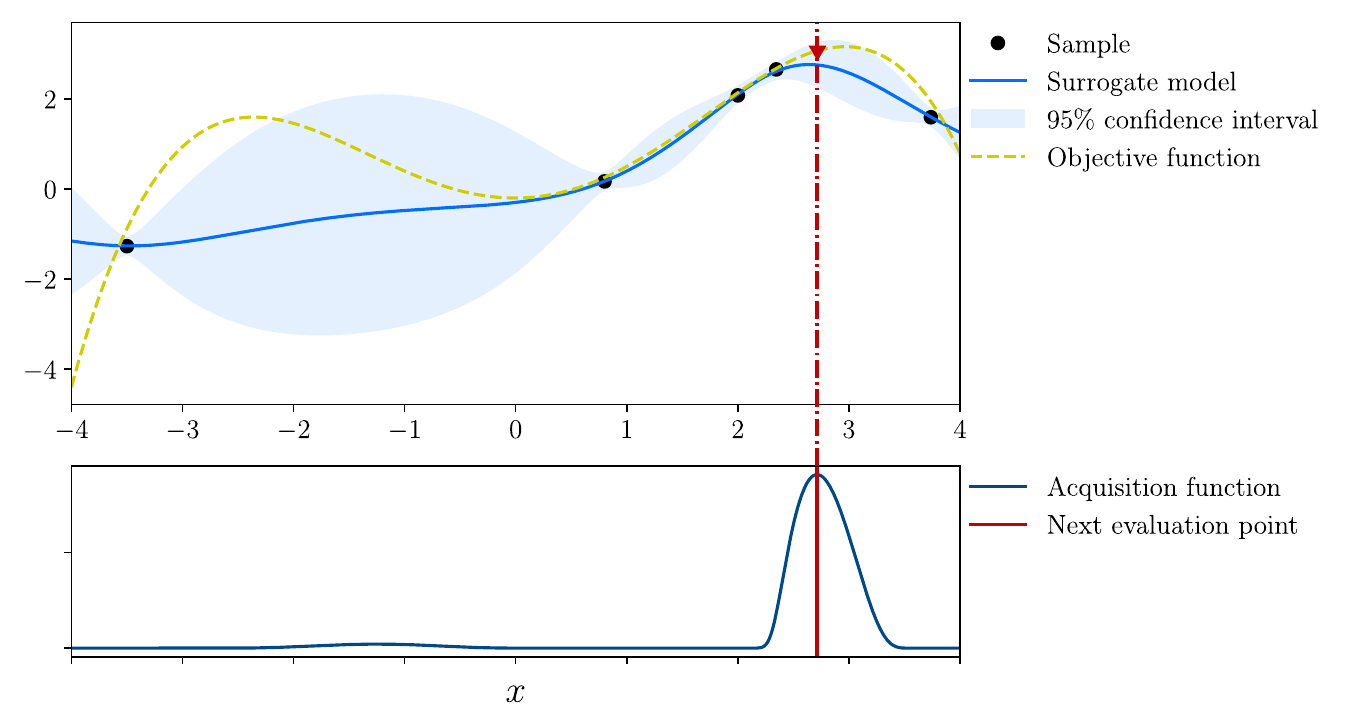}
	\caption{Illustration of the Bayesian optimisation process. The top figure shows the surrogate model (Gaussian process) fitted to the five data points from the objective function to be maximised. The bottom figure shows the acquisition function, where the red vertical line indicates the maximum value of the acquisition function. The red triangular marker is the data point to be added next. }
	\label{fig:BOAcquisition}
\end{figure}

In addition to the GP surrogate model, Bayesian optimisation requires an acquisition function. The Gaussian process regression starts with limited initial data. In each iteration during optimisation, a new data point $(\vec a, g)$ that maximises the acquisition function is determined and then added to the dataset, and the surrogate model is updated accordingly, as illustrated in Figure~\ref{fig:BOAcquisition}. In the presented examples, to improve the efficiency of GP regression, a sequential domain reduction scheme~\cite{Stander2002} provided in the Bayesian optimisation library~\cite{BayeOptPython2014} is employed.

Determining a new data point during optimisation entails a trade-off between exploration in the regions with high uncertainty and exploitation in the regions with high predicted objective function values.
The expected improvement (EI) acquisition function measures the expectation of improvement of the objective function at the next data point compared to the current best data. The maximum value of the objective function after $N$ iterations is denoted by $g(\vec{a}^*)$ with $\vec{a}^* = \underset{\vec{a}_i}{\mathop{\arg\max}} \, g(\vec{a}_i)$, where $i = 1, \cdots, N$. The EI acquisition function is derived analytically~\cite{Mockus1978, Jones1998} as
\begin{subequations}
	\begin{align}
		&EI(\vec{a}) = \delta(\vec{a}) P\left(\frac{\delta(\vec{a})}{\sigma(\vec{a})}\right) + \sigma(\vec{a}) p\left(\frac{\delta(\vec{a})}{\sigma(\vec{a})}\right) \, ,   \\
		&\delta(\vec{a}) = \mu(\vec{a}) - g(\vec{a}^*) - \xi \, ,
	\end{align}
\end{subequations}
where $p(\cdot)$ is the probability density function of the standard normal distribution and $P(\cdot)$ its cumulative density function, $\mu(\vec{a})$ and $\sigma(\vec{a})$ are the mean and standard deviation of the GP model evaluated at $\vec{a}$, $\delta(\vec{a})$ measures the improvement of $\mu(\vec{a})$ on $g(\vec{a}^*)$, and $\xi \in \mathbb R^+$ is a trade-off parameter that balances exploitation and exploration. A larger trade-off parameter $\xi$ may lead to more exploration, while a smaller $\xi$ may drive to more exploitation. In our implementation, the cross-sectional areas are scaled to the range between $0$ and $1$ when evaluating the acquisition function. Numerical experiments demonstrate that the value $\xi = 0.01$ gives good results~\cite{Lizotte2008}.


%% file: examples.tex
\section{Examples} \label{sec:examples}

\subsection{Two-ring star dome}
\begin{figure}
	\centering
	\subfloat[Top view] {
		\includegraphics[width=0.25 \textwidth]{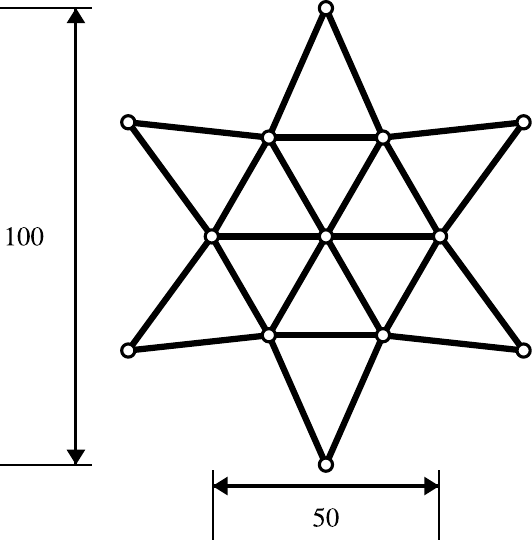}
		\label{fig:starDome_top} }
	\hfil
	\subfloat[Side view] {
		\includegraphics[width=0.32 \textwidth]{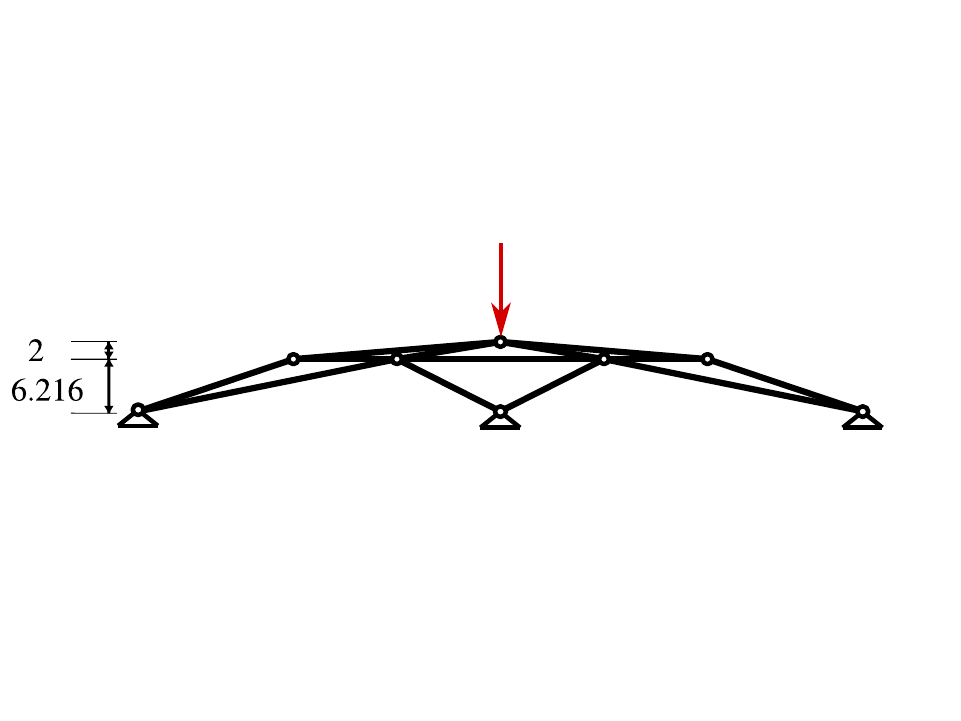}
		\label{fig:starDome_side} }
	\hfil	
	\subfloat[3D view] {
		\includegraphics[width=0.32 \textwidth]{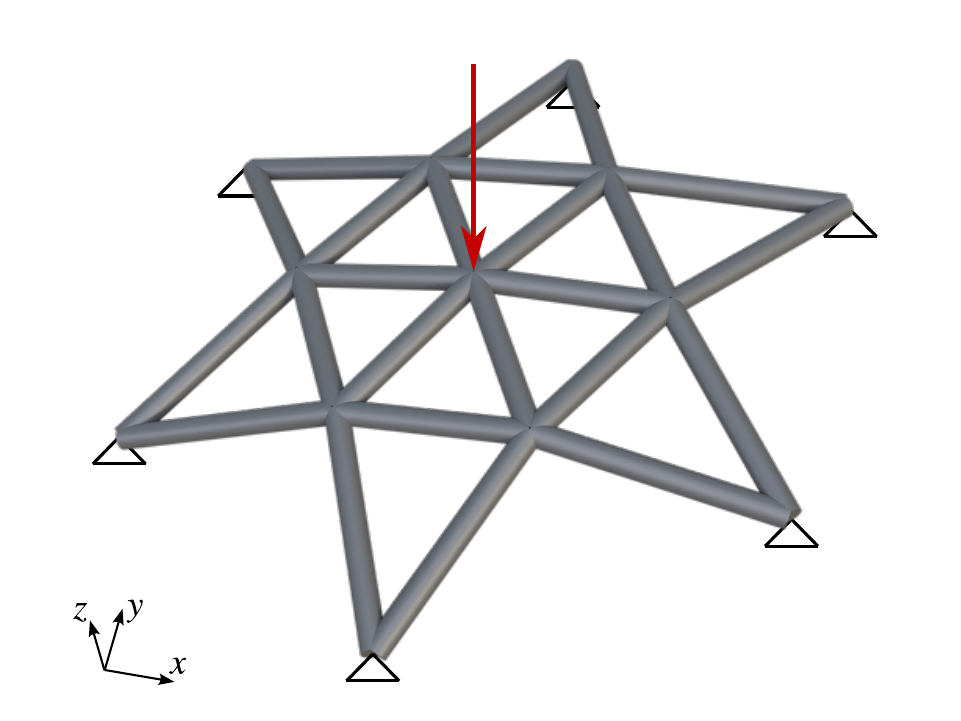}
		\label{fig:starDome_setup} }
	\caption{Geometry and boundary conditions of a two-ring star dome.}
	\label{fig:starDome}
\end{figure}

As the first example, we consider a 3D two-ring star dome shown in Figure~\ref{fig:starDome}. All the struts are assumed to be connected by pin joints. The pin supports are located at the six nodes of the outer ring, and a downward load is applied at the top middle node. The structure consists of $24$ elements and $13$ nodes. The Young's modulus of the material is $E=1.0 \times 10^8$ and the Poisson's ratio is $\nu=0.35$. The cross-sectional areas of the initial design are $a_i=0.5$ ($i = 1, 2, \cdots, 24$) for all struts, and the initial volume is $V_0=339.841$ which is considered as the volume constraint in optimisation.

\begin{figure}
	\centering
	\subfloat[Primary buckling mode shape] {
		\includegraphics[width=0.32 \textwidth]{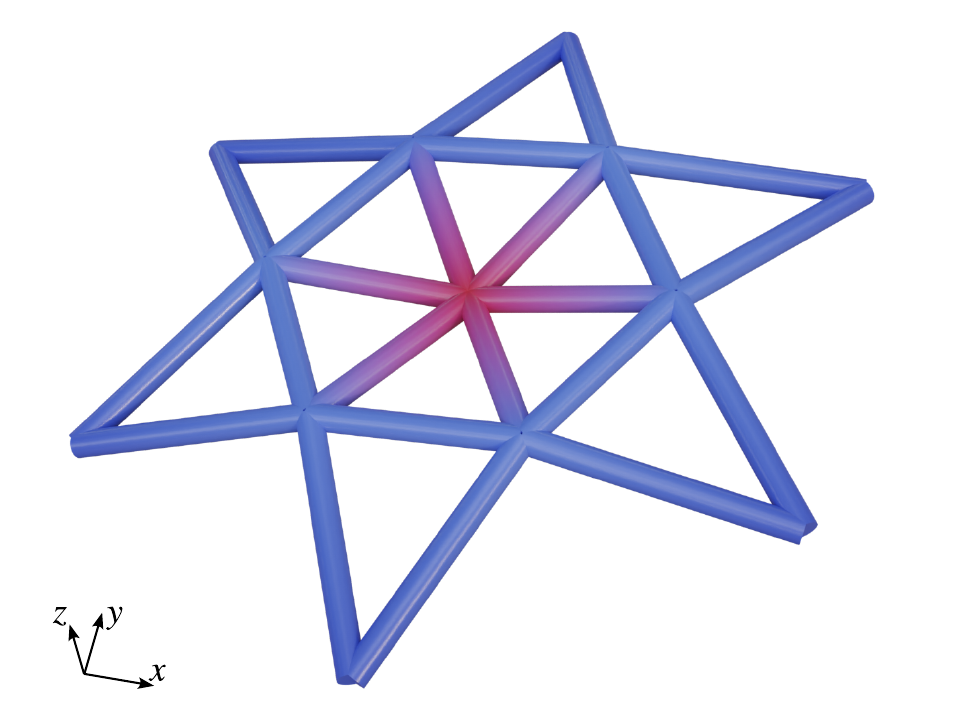}
		\label{fig:starDome_eigenModeShape} }
	\hfil
	\subfloat[Secondary buckling mode shape] {
		\includegraphics[width=0.32 \textwidth]{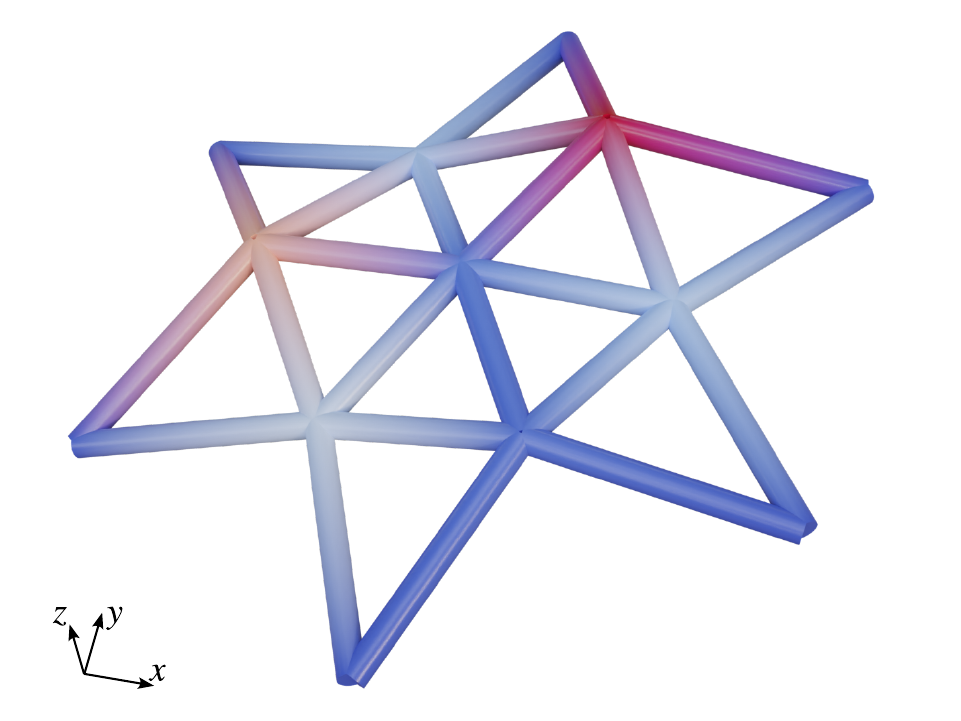}
		\label{fig:starDome_eigenModeShape2} }
	\hfil
	\subfloat[Design variables for optimisation] {
		\includegraphics[width=0.32 \textwidth]{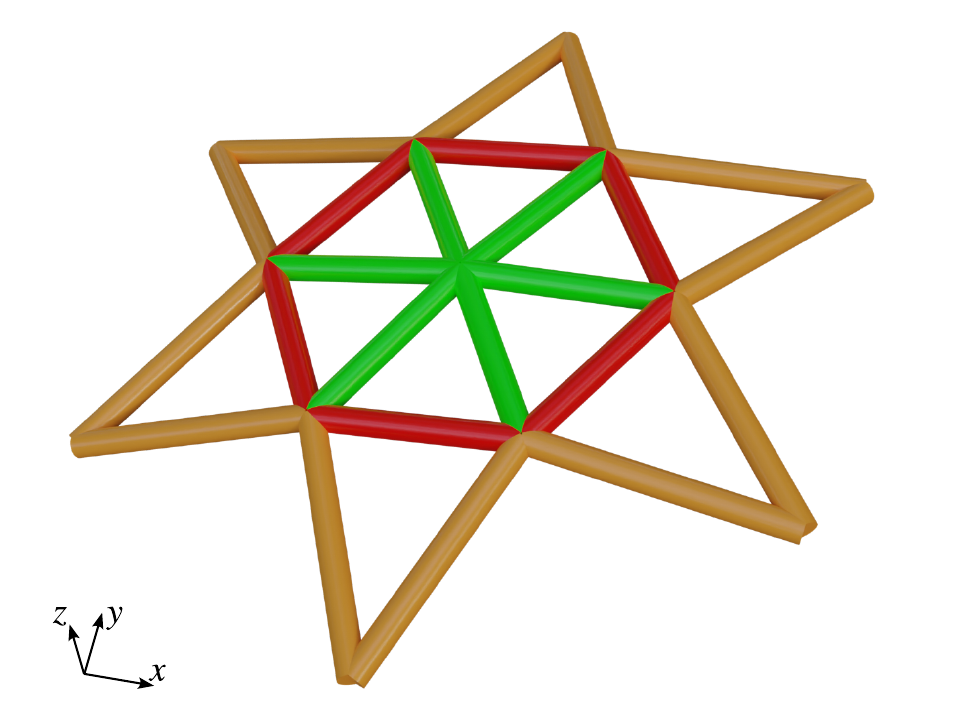}
		\label{fig:starDome_input_sameEle} }
	\caption{Mode shapes and design variables of the star dome. Figure (c) shows the three groups of cross-sectional areas in optimisation.}
	\label{fig:starDome_eigenModeShapeImperf}
\end{figure}

In robust optimisation the cross-sectional areas vary in the range $0.25 \leq a_i \leq 0.75$. Three scenarios of geometric imperfections are considered: primary buckling mode only ($\beta_1\vec{\phi}_1$), secondary buckling mode only ($\beta_2\vec{\phi}_2$) and the combination of the primary and secondary buckling modes ($\beta_1\vec{\phi}_1 + \beta_2\vec{\phi}_2$). The first two buckling modes are shown in Figures~\ref{fig:starDome_eigenModeShape} and~\ref{fig:starDome_eigenModeShape2}. The amplification factors $\beta_1$ and $\beta_2$ are sampled following a normal distribution $\mathcal{N}(0, 0.1^2)$. To reduce the dimensionality and computational cost of the optimisation problem, the struts are divided into three groups based on structural symmetry, such that the struts in the same group have identical cross-sectional areas, as shown in Figure~\ref{fig:starDome_input_sameEle} the three groups of struts are indicated in different colours. The maximum number of iterations in optimisation is set to $100$.

\begin{figure}
	\centering
	\subfloat[Pareto fronts and solution domain with random design variables] {
		\includegraphics[width=0.42 \textwidth]{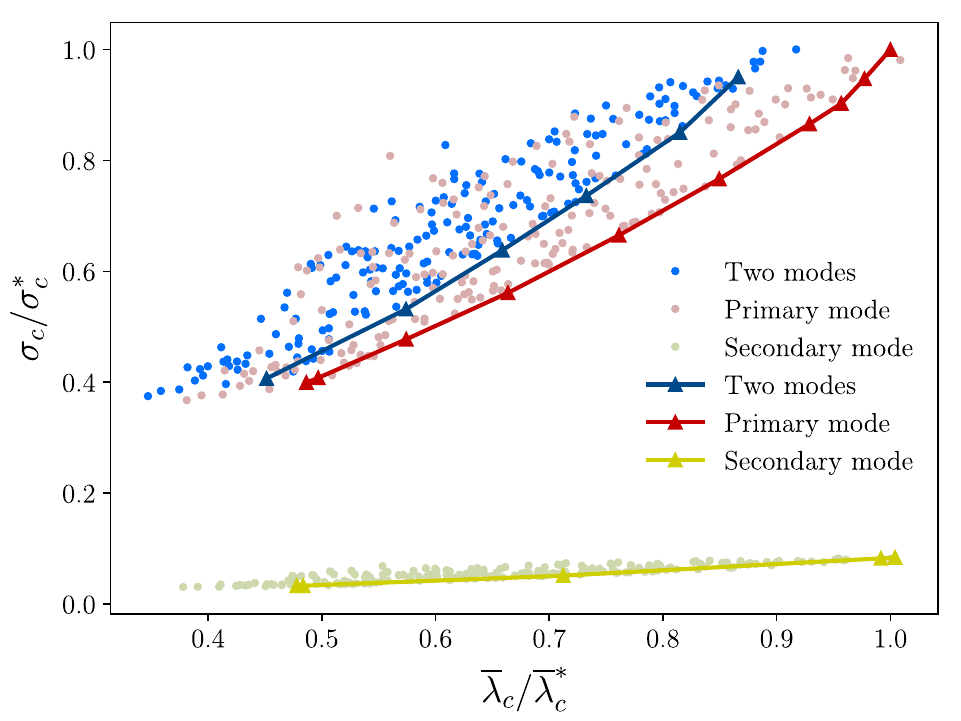}
		\label{fig:starDome_randomPareto} }
	\hfil
	\subfloat[Statistics of buckling loads] {
		\includegraphics[width=0.42 \textwidth]{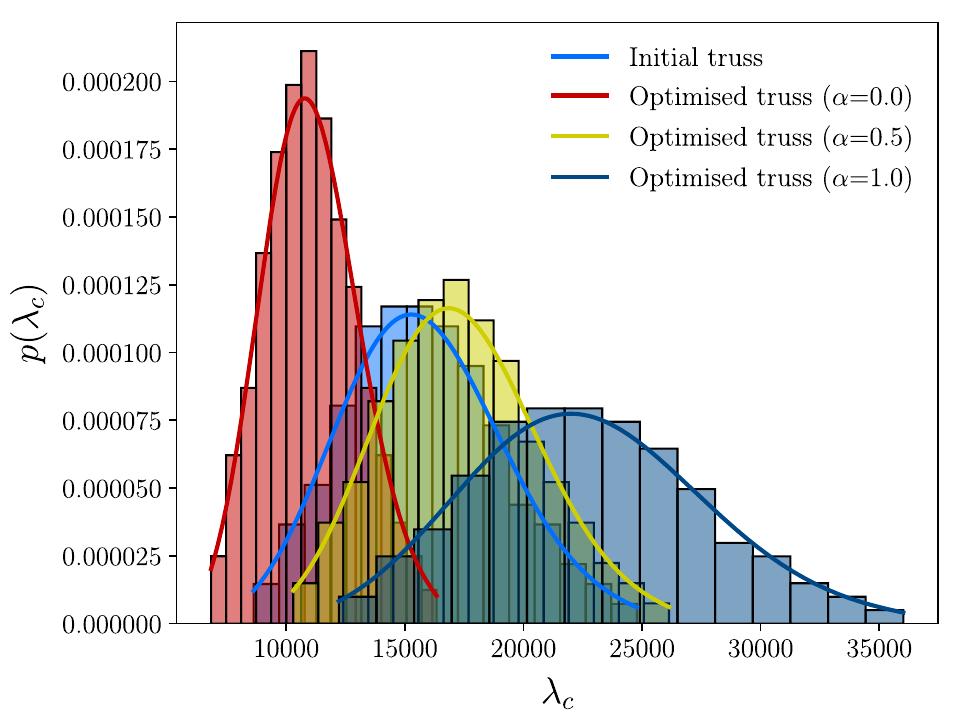}
		\label{fig:starDome_imperfDist_opt0.7+0.5} }
	\caption{Pareto fronts obtained with robust optimisation and buckling load statistics with different trade-off parameters $\alpha$. In (a) the solid lines represent the Pareto fronts, and the solution domain is sampled by evaluations with random design variables. In (b) the mean buckling loads are $15761$ (initial), $11055.78$ ($\alpha = 0$), $17319.41$ ($\alpha = 0.5$) and $22728.13$ ($\alpha = 1$); and the standard deviations are $3231.63$ (initial), $1901.02$ ($\alpha = 0$), $3166.89$ ($\alpha = 0.5$) and $4759.95$ ($\alpha = 1$).}
	\label{fig:starDome_pareto_imperfDist}
\end{figure}

\begin{figure*}
	\centering
	\subfloat[$\alpha=0.0$ ($\overline{\lambda}_c/\overline{\lambda}_c^* = 0.486$, $\sigma_{c}/\sigma_{c}^* = 0.399$)] {
		\includegraphics[width=0.25 \textwidth]{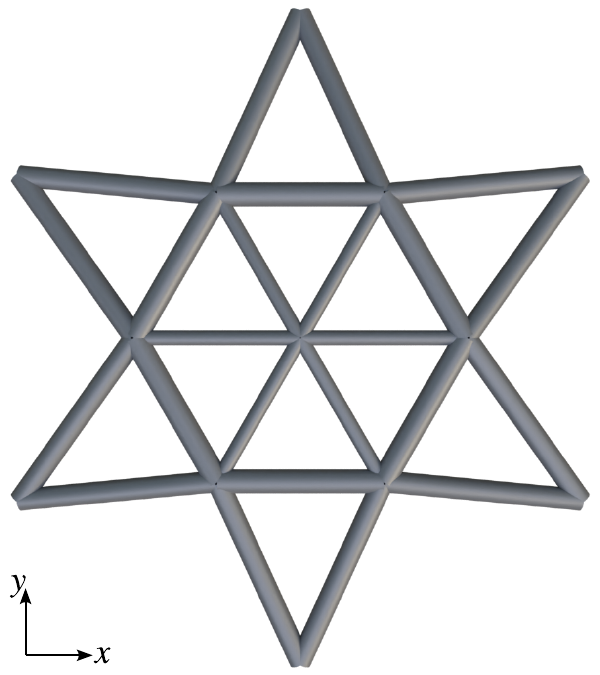}
		\label{fig:starDome_opt_alpha0.0} }
	\hfill
	\subfloat[$\alpha=0.5$ ($\overline{\lambda}_c/\overline{\lambda}_c^* = 0.761$, $\sigma_{c}/\sigma_{c}^* = 0.665$)] {
		\includegraphics[width=0.25 \textwidth]{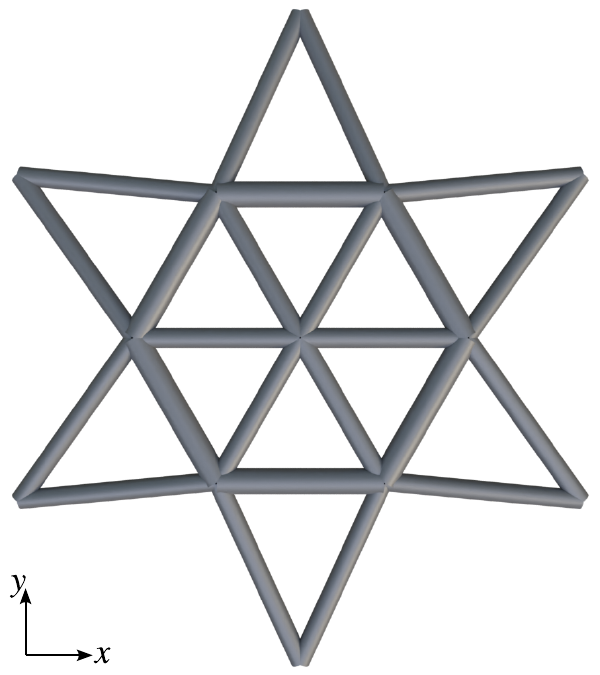}
		\label{fig:starDome_opt_alpha0.52} }
	\hfill
	\subfloat[$\alpha=1.0$ ($\overline{\lambda}_c/\overline{\lambda}_c^* = 1.0$, $\sigma_{c}/\sigma_{c}^* = 1.0$)] {
		\includegraphics[width=0.25 \textwidth]{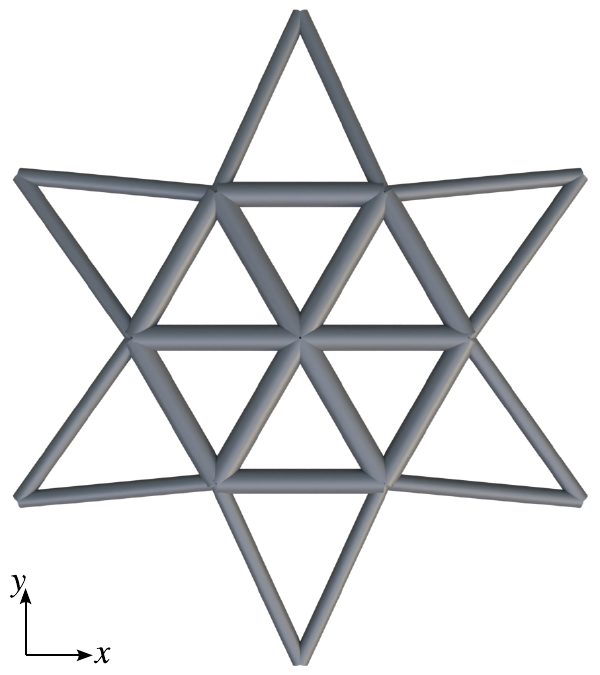}
		\label{fig:starDome_opt_alpha1.0} }
	\caption{Optimised cross-sectional areas of struts with different trade-off parameters $\alpha$.}
	\label{fig:starDome_opt_alpha}
\end{figure*}

Different values of the trade-off parameter $\alpha$ have been selected in robust optimisation, and a set of optimised results corresponding to the selected $\alpha$ form the Pareto front, as depicted in Figure~\ref{fig:starDome_randomPareto}. To verify the obtained Pareto front, a number of random design variables, i.e., the cross-sectional areas, are assigned to the struts, and the means and standard deviations of buckling loads corresponding to these assigned design variables form a solution domain, plotted as points in the figure. It turns out that the obtained Pareto front connecting optimised results follows closely the boundary of the solution domain. Comparing the Pareto fronts in the three scenarios of geometric imperfections, it is evident that geometric imperfections with the secondary buckling mode have much less effect than those with the primary buckling mode. To investigate the effect of $\alpha$, the buckling load statistics of the initial design and the three optimised results with $\alpha=0.0$, $\alpha=0.5$, $\alpha=1.0$ are compared in Figure~\ref{fig:starDome_imperfDist_opt0.7+0.5}, in which the primary buckling mode is considered for geometric imperfections. For different $\alpha$ values, the resulting buckling loads have distinct mean and standard deviation values. For $\alpha=0.0$, the robust optimisation problem minimises the standard deviation only; as $\alpha$ increases, both the mean and standard deviation of buckling loads increase; for $\alpha=1.0$, it maximises the mean value only, equivalent to the deterministic optimisation for buckling load; on the other hand, for $\alpha = 0.5$, the mean buckling load is larger than that of the initial truss, and the standard deviation is smaller than that of the initial truss, yielding overall a better structure. The number of objective function evaluations in Bayesian optimisation to reach the optimum is 42 in the case of $\alpha = 0.5$. The optimised struts of the dome considering $\alpha=0.0$, $0.5$ and $1.0$ are shown in Figure~\ref{fig:starDome_opt_alpha}. The cross-sectional areas of struts in the inner ring increase along with $\alpha$, whereas those in the outer ring decrease.

\subsection{Five-ring star dome}

\begin{figure*}
	\centering
	\subfloat[Top view] {
		\includegraphics[width=0.28 \textwidth]{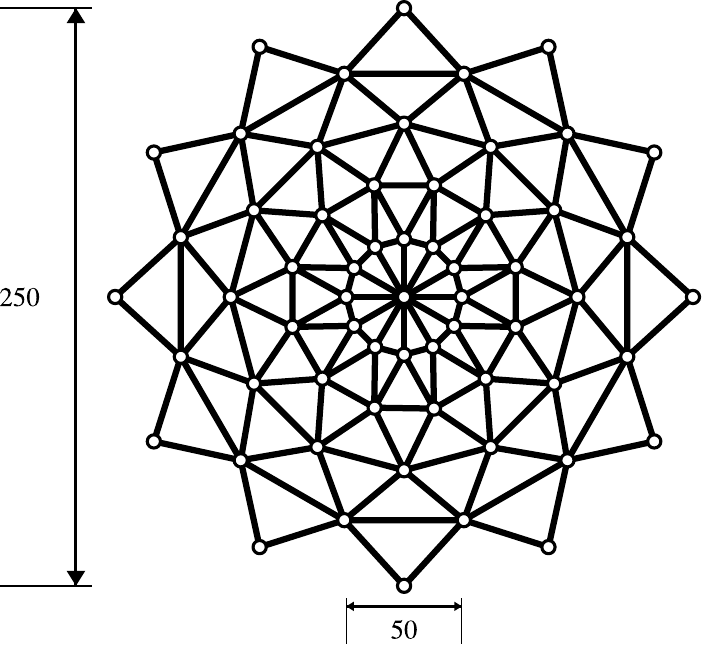}
		\label{fig:starDome5_top} }
	\hfil
	\subfloat[Side view] {
		\includegraphics[width=0.32 \textwidth]{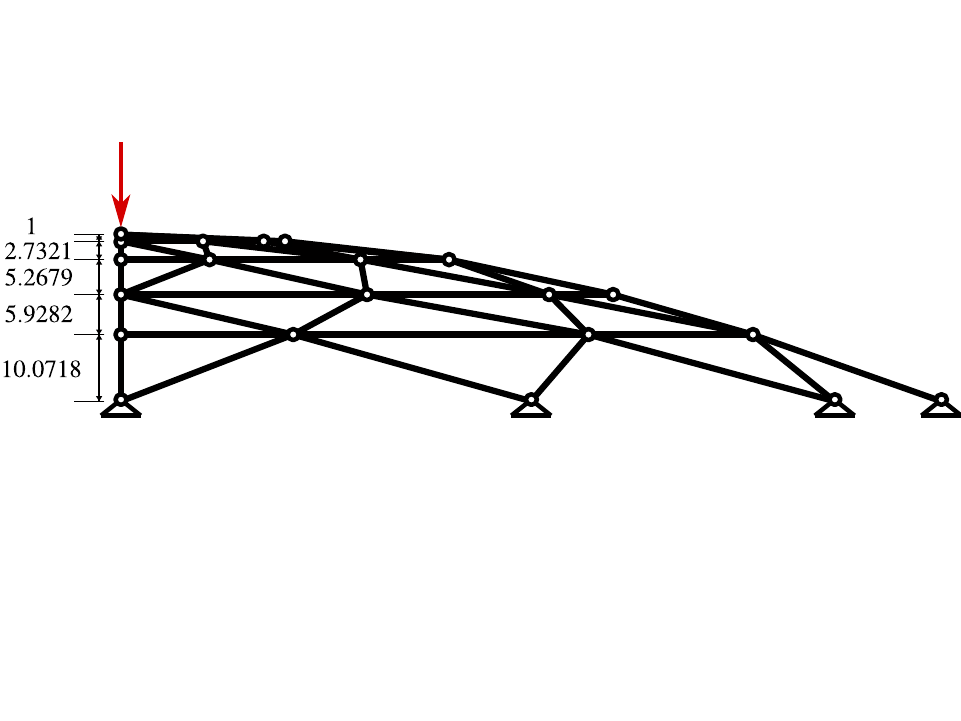}
		\label{fig:starDome5_side} }
	\hfil	
	\subfloat[3D view] {
		\includegraphics[width=0.32 \textwidth]{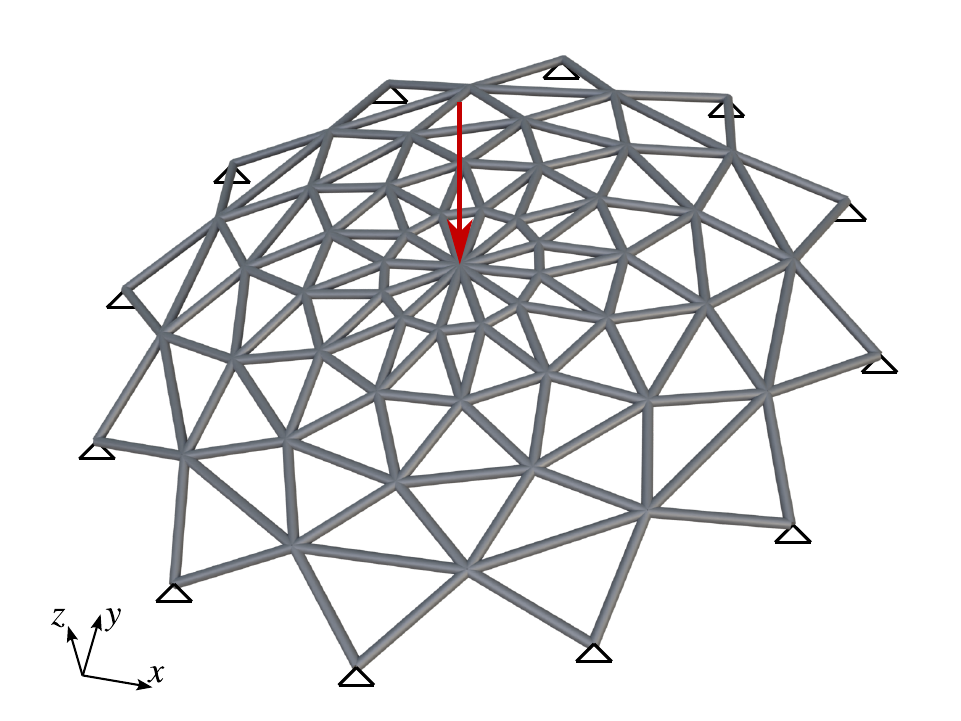}
		\label{fig:starDome5_setup} }
	\caption{Geometry and boundary conditions of a five-ring star dome.}
	\label{fig:starDome5}
\end{figure*}

Next, a five-ring star dome with all the struts pin jointed is considered, as shown in Figure~\ref{fig:starDome5}. The structure consists of $156$ elements and $61$ nodes. The material properties of all struts remain the same as the two-ring star dome. The cross-sectional areas of all struts in the initial design are $a_i=0.5$ ($i = 1, 2, \cdots, 156$). The initial volume of the dome is $V_0=2497.12$ and is considered as the volume constraint in optimisation.

\begin{figure*}
	\centering
	\subfloat[Primary buckling mode shape] {
		\includegraphics[width=0.32 \textwidth]{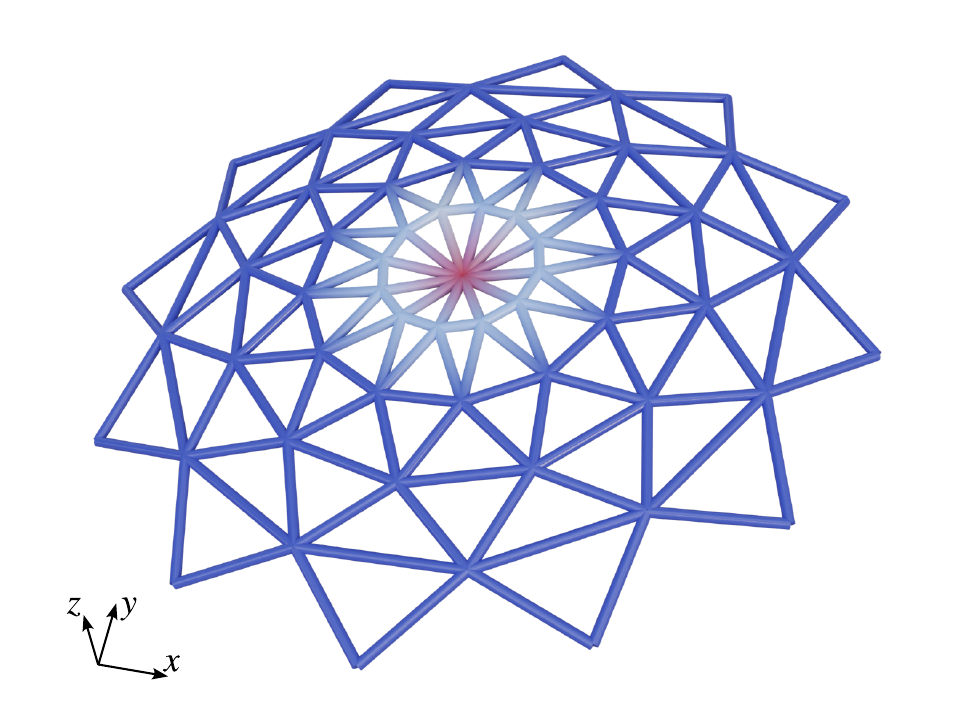}
		\label{fig:starDome5_eigenModeShape} }
	\hfil
	\subfloat[Secondary buckling mode shape] {
		\includegraphics[width=0.32 \textwidth]{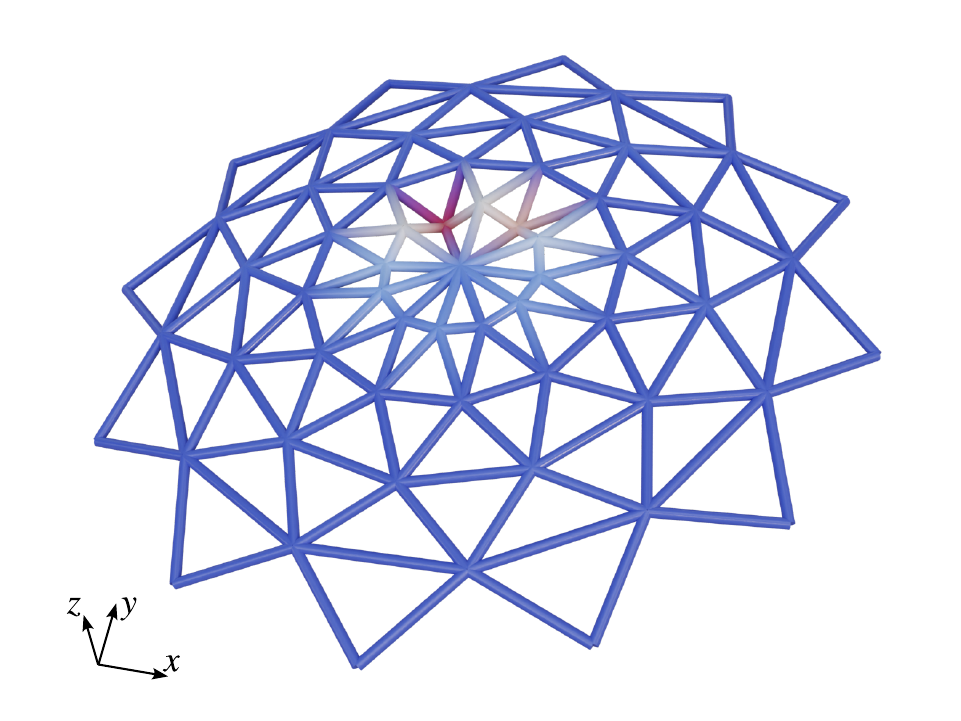}
		\label{fig:starDome5_eigenModeShape2} }
	\hfil
	\subfloat[Design variables for optimisation] {
		\includegraphics[width=0.32 \textwidth]{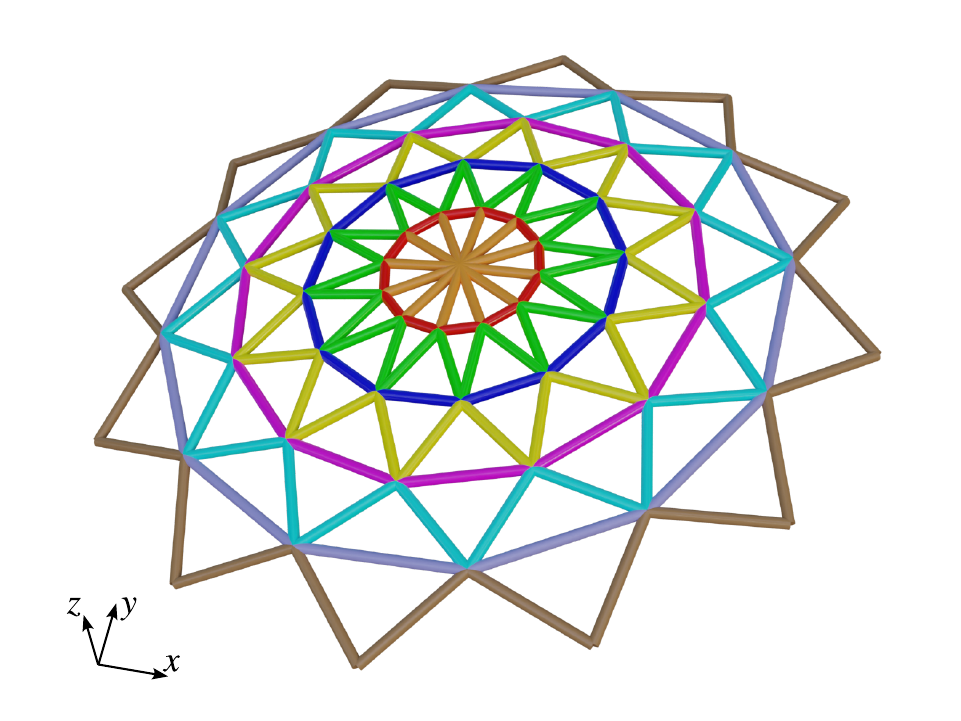}
		\label{fig:starDome5_input_sameEle} }
	\caption{Mode shapes and design variables of the star dome. Figure (c) shows the nine groups of cross-sectional areas in optimisation.}
	\label{fig:starDome5_eigenModeShapeImperf}
\end{figure*}

The same as the two-ring star dome, the cross-sectional areas of struts vary in the range $0.25 \leq a_i \leq 0.75$ in robust optimisation, and two scenarios of geometric imperfections are considered: primary buckling mode only ($\beta_1\vec{\phi}_1$) and secondary buckling mode only ($\beta_2\vec{\phi}_2$). The first two buckling modes of the initial design are shown in Figure~\ref{fig:starDome5_eigenModeShape} and~\ref{fig:starDome5_eigenModeShape2}. The amplification factors $\beta_1$ and $\beta_2$ are sampled following a normal distribution $\mathcal{N}(0, 0.03^2)$. Based on structural symmetry, the cross-sectional areas of struts are divided into nine groups to have nine design variables in optimisation, as depicted in Figure~\ref{fig:starDome5_input_sameEle}. The maximum number of optimisation iterations is set to $200$. Compared with the optimisation of the two-ring star dome, a much larger number of degrees-of-freedom and design variables are involved in the five-ring star dome, resulting in more expensive objective function evaluations and more optimisation steps.

\begin{figure*}
	\centering
	\subfloat[Random samples and Pareto fronts] {
		\includegraphics[width=0.42 \textwidth]{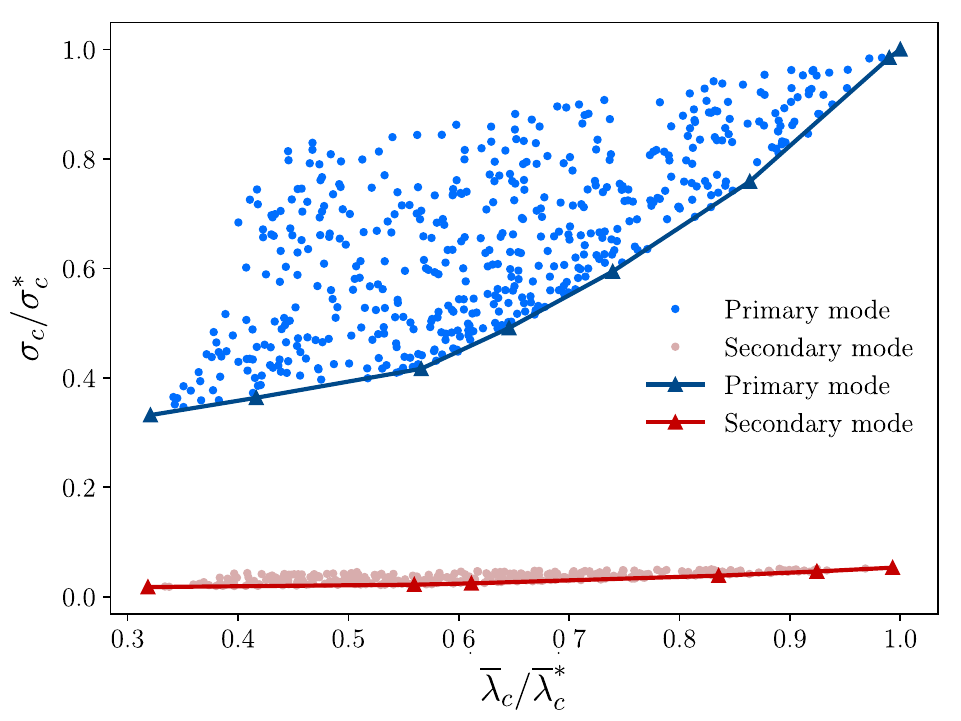}
		\label{fig:starDome5_randomParetoFront} }
	\hfil
	\subfloat[Selected buckling load distributions] {
		\includegraphics[width=0.42 \textwidth]{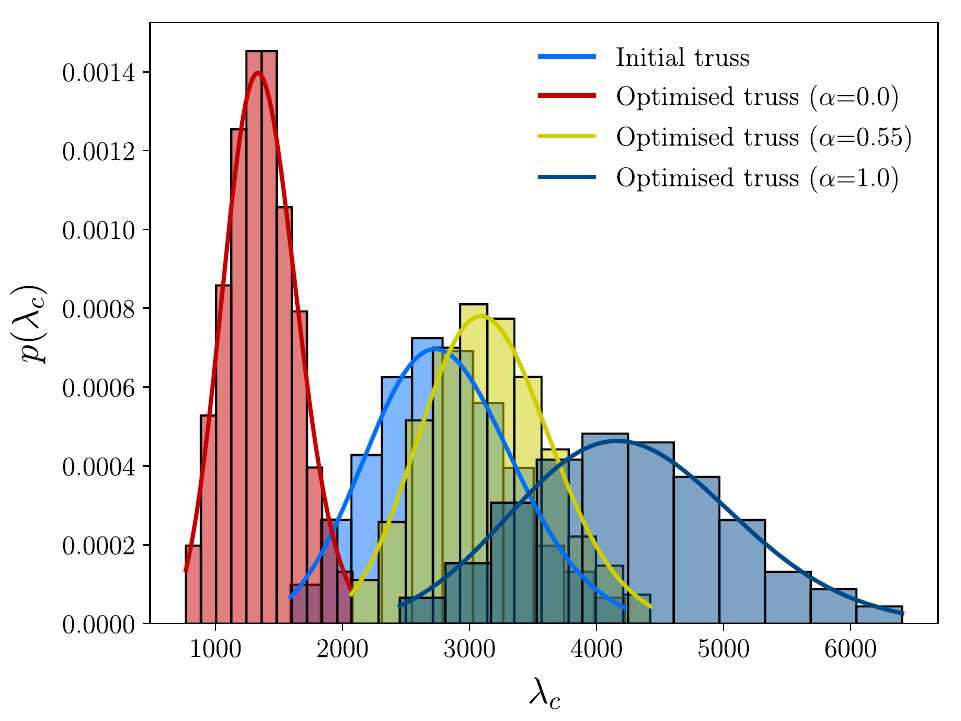}
		\label{fig:starDome5_imperfDist_opt0.6+0.55} }
	\caption{Pareto fronts obtained with robust optimisation and buckling load statistics with different trade-off parameters $\alpha$. In (a) the solid lines represent the Pareto fronts, and the solution domain is sampled by evaluations with random design variables. In (b) the mean buckling loads are $2790.56$ (initial), $1363.96$ ($\alpha = 0$), $3144.95$ ($\alpha = 0.55$) and $4253.12$ ($\alpha = 1$); and the standard deviations are $526.66$ (initial), $262.62$ ($\alpha = 0$), $470.29$ ($\alpha = 0.55$) and $791.40$ ($\alpha = 1$).}
	\label{fig:starDome5_pareto_imperfDist}
\end{figure*}

\begin{figure*}[h!]
	\centering
	\subfloat[$\alpha=0.0$ ($\overline{\lambda}_c/\overline{\lambda}_c^* = 0.321$, $\sigma_{c}/\sigma_{c}^* = 0.332$)] {
		\includegraphics[width=0.3 \textwidth]{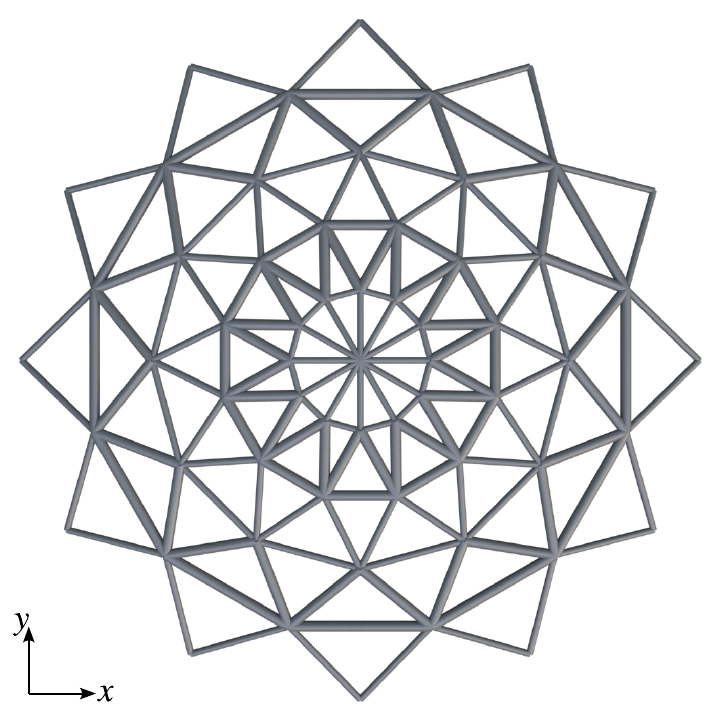}
		\label{fig:starDome5_opt_render_alpha0.0} }
	\hfill
	\subfloat[$\alpha=0.55$ ($\overline{\lambda}_c/\overline{\lambda}_c^* = 0.739$, $\sigma_{c}/\sigma_{c}^* = 0.594$)] {
		\includegraphics[width=0.3 \textwidth]{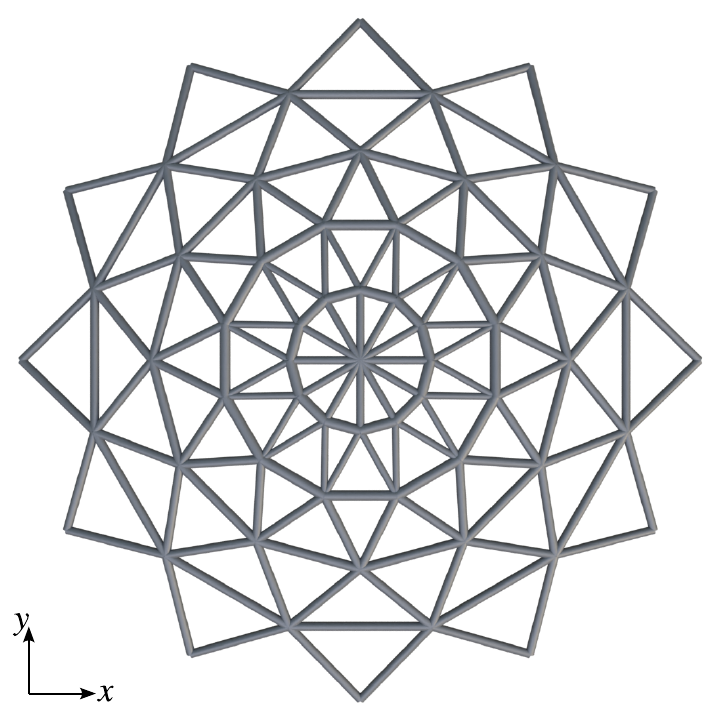}%
		\label{fig:starDome5_opt_render_alpha0.5} }
	\hfill
	\subfloat[$\alpha=1.0$ ($\overline{\lambda}_c/\overline{\lambda}_c^* = 1.0$, $\sigma_{c}/\sigma_{c}^* = 1.0$)] {
		\includegraphics[width=0.3 \textwidth]{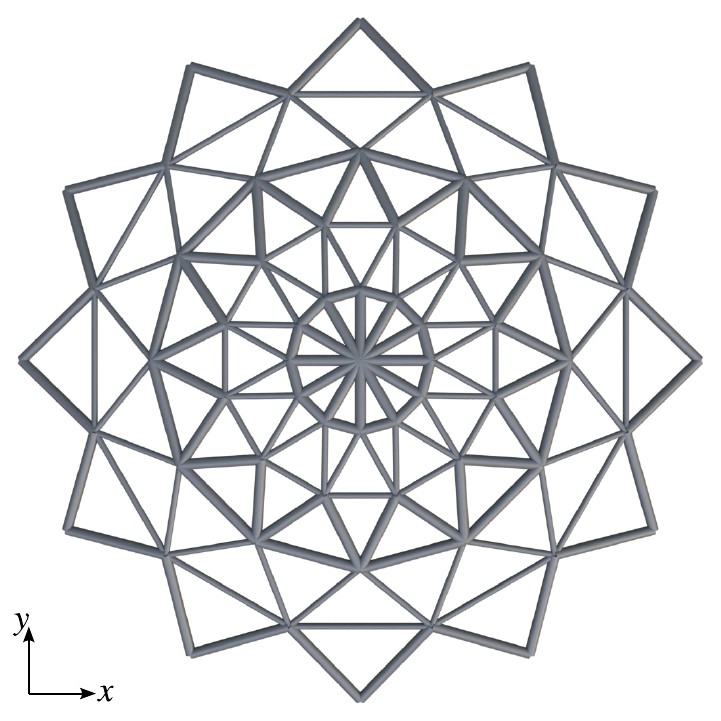}%
		\label{fig:starDome5_opt_render_alpha1.0} }
	\caption{Optimised cross-sectional areas of struts with different trade-off parameters $\alpha$.}
	\label{fig:starDome5_opt_render_alphaCompare}
\end{figure*}

The Pareto fronts of the optimised results considering different values of $\alpha$ are displayed in Figure~\ref{fig:starDome5_randomParetoFront}. Again, a number of random design variables are considered to compute the means and standard deviations of buckling loads to verify the obtained Pareto fronts. The comparison of the buckling load statistics of the initial design and the optimised structures with $\alpha=0.0$, $0.55$ and $1.0$ are plotted in Figure~\ref{fig:starDome5_imperfDist_opt0.6+0.55}, in which the primary buckling mode is considered for geometric imperfections. The corresponding optimised cross-sectional areas of struts are shown in Figure~\ref{fig:starDome5_opt_render_alphaCompare}. It can be seen that the cross-sectional areas of struts are more uniform with $\alpha=0.55$ than those with $\alpha=0.0$ and $1.0$. In the case of $\alpha = 0.55$, the buckling load has a larger mean and a smaller standard deviation than those of the initial truss, yielding overall a better structure. The number of objective function evaluations in Bayesian optimisation to reach the optimum is 54 in this case.

\subsection{Truss column}
\begin{figure}[h]
\begin{minipage}{0.48 \textwidth}
	\centering
	\subfloat[Top view] {
		\includegraphics[width=0.25 \textwidth]{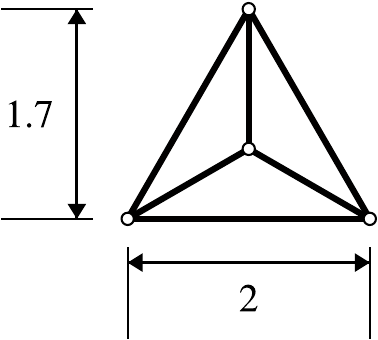}
		\label{fig:trussColumn_top} }
	\hfil
	\subfloat[Side view] {
		\includegraphics[width=0.215 \textwidth]{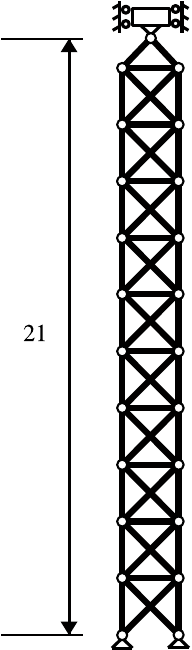}
		\label{fig:trussColumn_side} }
	\hfil	
	\subfloat[3D view] {
		\includegraphics[width=0.4 \textwidth]{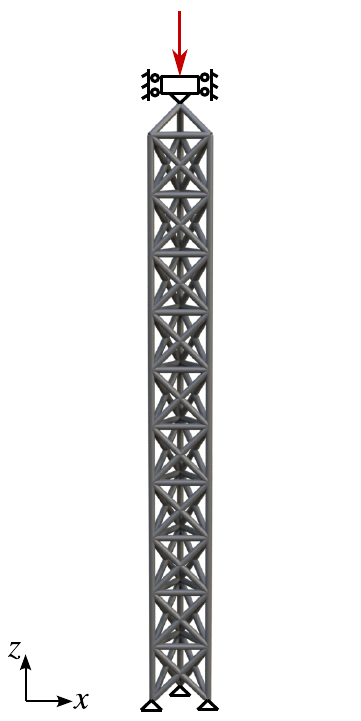}
		\label{fig:trussColumn_setup} }
	\caption{Geometry and boundary conditions of a truss column in different views.}
	\label{fig:trussColumn}
\end{minipage}
\hfill
\begin{minipage}{0.45 \textwidth}
	\centering
	\subfloat[Buckling mode shape] {
		\includegraphics[width=0.455 \textwidth]{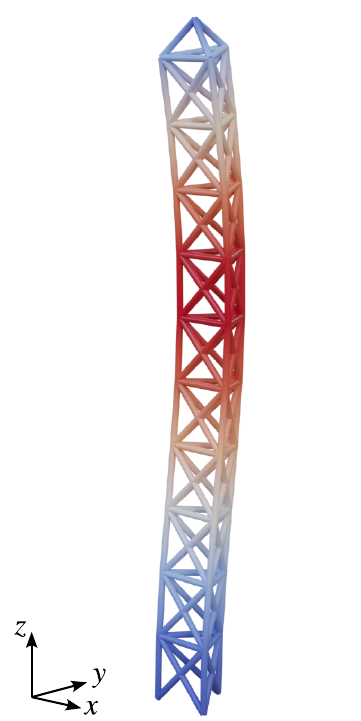}
		\label{fig:trussColumn_eigenModeShape} }
	\hfil
	\subfloat[Design variables for optimisation] {
		\includegraphics[width=0.455 \textwidth]{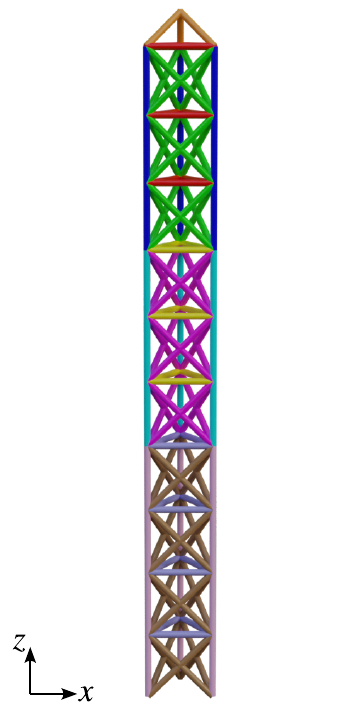}
		\label{fig:trussColumn_input_sameEle} }
	\caption{Mode shapes and design variables of the truss column. Figure (b) shows the ten groups of cross-sectional areas in optimisation.}
	\label{fig:trussColumn_imperf_bucklingF}
\end{minipage}
\end{figure}

A truss column structure is further considered, with ten blocks of triangular prisms and a tetrahedron on the top, as shown in Figure~\ref{fig:trussColumn}. It is pinned to the ground at three bottom nodes. The top node is constrained horizontally with a downward load. The structure has $123$ elements and $34$ nodes. The Young's modulus of the material is $E=1.0\times10^4$ and Poisson's ratio is $\nu=0.35$. The cross-sectional areas of the initial design are $a_i=0.1$ ($i = 1, 2, \cdots, 123$) for all struts. The initial volume is $V_0=29.3$ which is the volume constraint in optimisation.

In robust optimisation, the cross-sectional areas vary in the range $0.05 \leq a_i \leq 0.15$. The primary buckling mode (see Figure~\ref{fig:trussColumn_eigenModeShape}) is considered to generate geometric imperfections ($\beta_1\vec{\phi}_1$) with the amplification factor $\beta_1$ following a normal distribution $\mathcal{N}(0, 0.006^2)$. The cross-sectional areas of struts are divided into $10$ groups with each representing one design variable, as depicted in Figure~\ref{fig:trussColumn_input_sameEle}. In the lower four blocks, the cross-sectional areas of the vertical struts, horizontal struts and cross-bracing struts are in three groups. Six other groups are set to the fifth to seventh blocks and eighth to tenth blocks. The last group includes the three struts of the tetrahedron on the top. The maximum number of iterations in optimisation is set to $300$.

\begin{figure*}
	\centering
	\subfloat[Random samples and Pareto front] {
		\includegraphics[width=0.45 \textwidth]{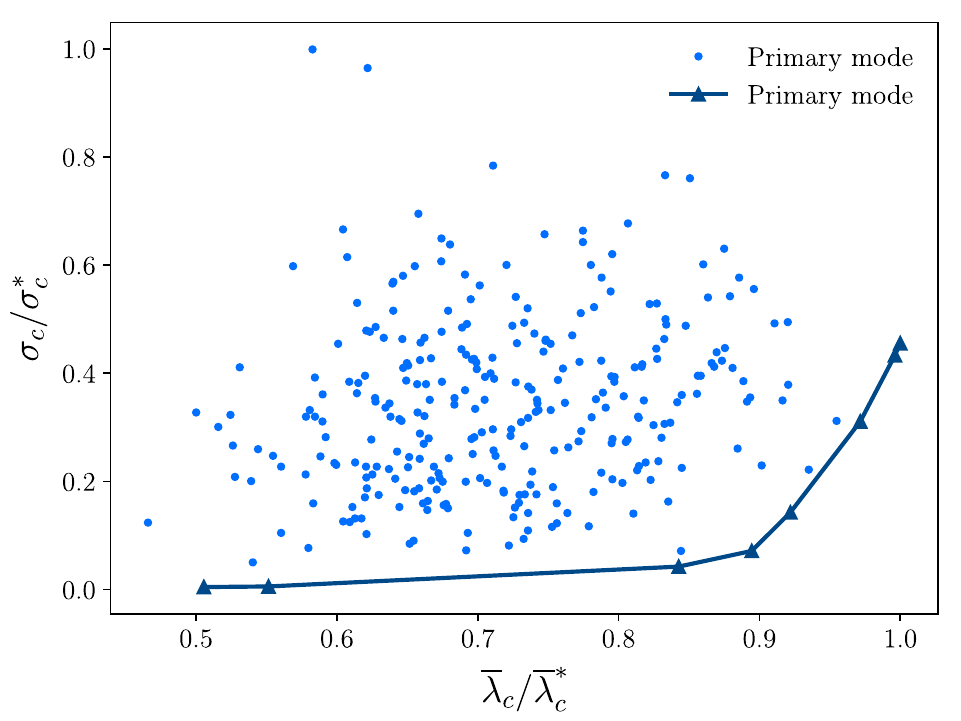}
		\label{fig:trussColumn_randomParetoFront} }
	\hfil
	\subfloat[Selected buckling load distributions] {
		\includegraphics[width=0.45 \textwidth]{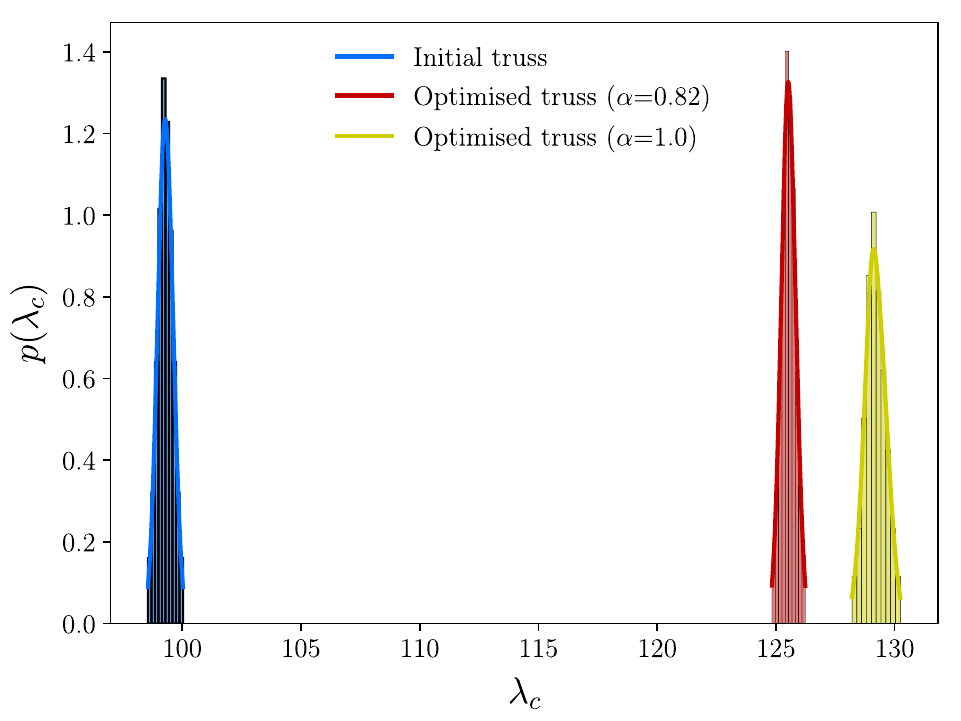}
		\label{fig:trussColumn_imperfDist_opt0.82+0.85} }
	\caption{Pareto fronts obtained with robust optimisation and buckling load statistics with different trade-off parameters $\alpha$. In (a) the solid lines represent the Pareto fronts, and the solution domain is sampled by evaluations with random design variables. In (b) the mean buckling loads are $99.29$ (initial), $125.53$ ($\alpha = 0.82$) and $129.2$ ($\alpha = 1$); and the standard deviations are $0.295$ (initial), $0.279$ ($\alpha = 0.82$) and $0.405$ ($\alpha = 1$).}
	\label{fig:trussColumn_pareto_imperfDist}
\end{figure*}

\begin{figure}
	\centering
	\subfloat[$\alpha=0.0$ ($\overline{\lambda}_c/\overline{\lambda}_c^* = 0.506$, $\sigma_{c}/\sigma_{c}^* = 0.004$)] {
		\includegraphics[width=0.25 \textwidth]{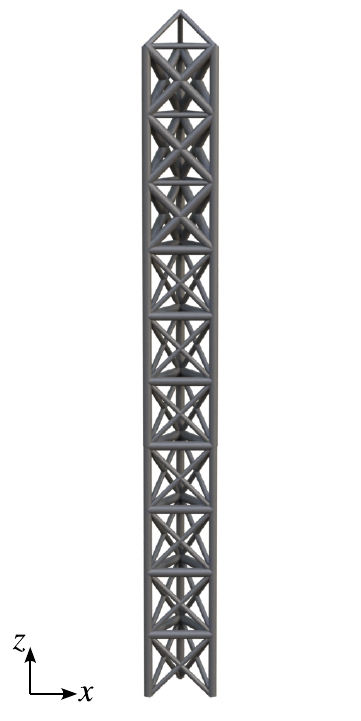}
		\label{fig:trussColumn_opt_render_alpha0.0} }
	\hfil
	\subfloat[$\alpha=0.82$ ($\overline{\lambda}_c/\overline{\lambda}_c^* = 0.972$, $\sigma_{c}/\sigma_{c}^* = 0.311$)] {
		\includegraphics[width=0.25 \textwidth]{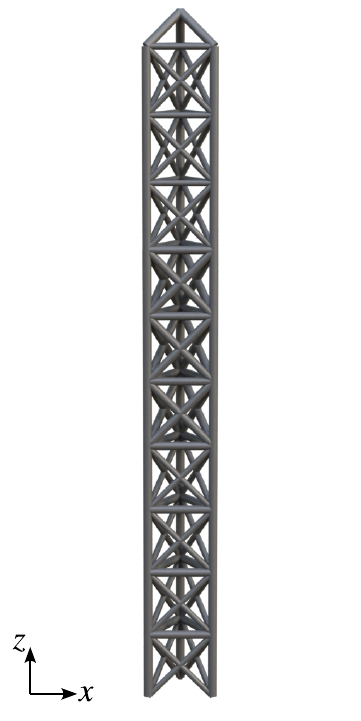}%
		\label{fig:trussColumn_opt_render_alpha0.65} }
	\hfil
	\subfloat[$\alpha=1.0$ ($\overline{\lambda}_c/\overline{\lambda}_c^* = 1.0$, $\sigma_{c}/\sigma_{c}^* = 0.455$)] {
		\includegraphics[width=0.25 \textwidth]{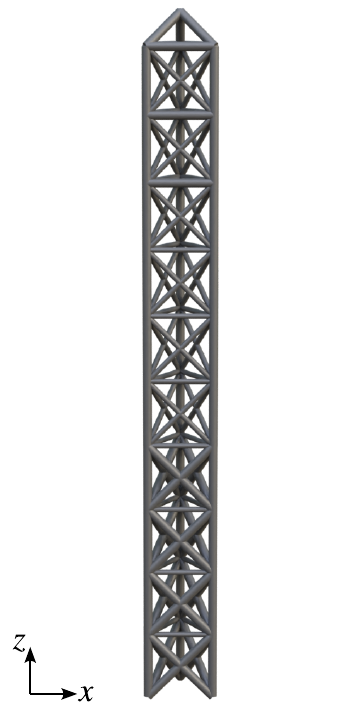}%
		\label{fig:trussColumn_opt_render_alpha1.0} }
	\caption{Optimised cross-sectional areas of struts with different trade-off parameters $\alpha$.}
	\label{fig:trussColumn_opt_render_alphaCompare}
\end{figure}

The Pareto front of the optimised results is shown in Figure~\ref{fig:trussColumn_randomParetoFront}. The buckling load statistics of the initial design and the optimised structures with $\alpha=0.82$ and $1.0$ are shown in Figure~\ref{fig:trussColumn_imperfDist_opt0.82+0.85}. For $\alpha=0.0$, the buckling load has a much smaller mean and standard deviation, which is almost a vertical line on the left of the plot and hence not included in the plot. The three optimised truss columns with different $\alpha$ values are compared in Figure~\ref{fig:trussColumn_opt_render_alphaCompare}. In the case of $\alpha=0.0$, thicker struts are in the top blocks and thinner struts in the middle and bottom blocks; in the case of $\alpha=1.0$, thicker bars are in the bottom blocks; in the case of $\alpha=0.82$, thicker bars are in the middle blocks.

In this example, the effectiveness of the proposed robust optimisation for the problem up to $10$ dimensions has been validated. For $\alpha = 0.82$, the buckling load of the truss has a larger mean and a smaller standard deviation than the initial truss, yielding overall a superior structure. In this case, the number of objective evaluations in Bayesian optimisation to reach the optimum is 75. Different from the star domes which exhibit snap-through buckling, the truss column exhibits bifurcation buckling. Hence, the proposed robust optimisation is valid for different types of stability points.

%% file: conclusions.tex
\section{Conclusions \label{sec:conclusions}}
%
We introduced a novel framework for robust optimisation of buckling loads of structures with geometric uncertainties resulting from random geometric imperfections and demonstrated its application to the sizing optimisation of spatial trusses.  The geometric imperfections are parameterised in terms of the random amplitudes of the as-designed structure's buckling modes. The prescribed Gaussian imperfection parameters, which are statistically independent, have a prescribed mean and standard deviation. The corresponding non-Gaussian buckling load distribution is determined by first sampling from the imperfections and then performing a geometrically nonlinear finite element analysis for each sample. The buckling loads are computed with an efficient extended system method, avoiding costly path-following procedures directly. For sampling from the imperfection parameter distribution, we employ the quasi-Monte Carlo sampling, i.e. Sobol sampling, which uses a quasi-random sequence to provide samples with a greater spatial uniformity than purely random sampling.  After obtaining the samples, a Gaussian distribution is obtained by transforming the samples using the Gaussian cumulative distribution function.   Since the generated samples are more uniformly distributed than random samples, the number of samples needed to characterise the probability distribution and the number of finite element evaluations are both significantly reduced. Furthermore, we use the gradient-free Bayesian optimisation, which is exceptionally suitable given that the objective function is expensive and the number of design variables is small.

We first verified the effectiveness and accuracy of the proposed algorithm for determining the buckling load distribution using a classic Von Mises truss example with analytical solutions. After that, we demonstrated the efficiency and efficacy of our robust optimisation framework with three spatial truss examples having an increasing number of design variables.   Additionally, we studied the effect of the trade-off parameter on robust optimisation, comparing the results for different trade-off parameter values, and determining the Pareto fronts. In our numerical examples, the cross-sectional areas of struts are grouped to a small number of design variables to balance the structural performance and computational efficiency. Since Gaussian process is used to model the objective function in Bayesian optimisation and the statistics of buckling loads are computed with quasi-Monte Carlo sampling, the scalability of the approach to problems of higher dimensions requires additional considerations. Practical techniques to cope with higher dimensional problems include dimensionality-reduction techniques~\cite{archbold2024surrogate,vadeboncoeur2023fully} and batch-based parallelism~\cite{Wang2018} to allow the large-scale Bayesian optimisation; see also Wang et al.~\cite{Wang2023} for a recent review on advanced techniques in Bayesian optimisation.

In closing, we note some promising directions for future research and extensions of the proposed approach. Although only the sizing optimisation of trusses was considered, the framework can also be applied to other kinds of structures, such as beams, shells and frames, and extended to shape or topology optimisation, see e.g.~\cite{bandara2018isogeometric, Xiao2019, XiaoCirak2022}. Similarly, it is straightforward to consider other than geometric uncertainties, such as imperfections in material properties and loading. Although the presented results were obtained with an in-house nonlinear finite element code, the modular structure makes the approach appealing for implementation as a user-defined subroutine in commercial software packages, like Abaqus or Ansys.  Furthermore, with sufficient training data and computational resources, one could consider employing a pre-trained neural network or Gaussian process surrogate to predict the entire nonlinear structural behaviour~\cite{archbold2024surrogate,vadeboncoeur2023fully}.